%% file: main.tex
\documentclass[a4paper,11pt,hidelinks]{article}

\usepackage{subfiles}
\newcommand{\onlyinsubfile}[1]{#1}
\newcommand{\notinsubfile}[1]{}

\usepackage{amsmath,multicol}
\usepackage{amssymb,mathrsfs}
\usepackage{amsthm}
\usepackage[all]{xypic}
\usepackage{amsfonts}
\usepackage{fullpage}
\usepackage{verbatim}
\usepackage{graphicx}
\usepackage{enumitem}
\usepackage{thmtools}
\usepackage{url}  
\usepackage[dvipsnames]{xcolor}
\usepackage{faktor}
\usepackage{tikz-cd}
\usepackage{ifthen}

\usepackage{hyperref}
\usepackage[nameinlink]{cleveref}
\crefname{subsection}{subsection}{subsections}

\numberwithin{equation}{section}

\theoremstyle{plain}

\newtheorem{thm}[equation]{Theorem}
\newtheorem{proposition}[equation]{Proposition}
\newtheorem{corollary}[equation]{Corollary}
\newtheorem{lemma}[equation]{Lemma}
\newtheorem{question}[equation]{Question}

\newtheorem{claim}[equation]{Claim}

\theoremstyle{definition}
\newtheorem{defn}[equation]{Definition}

\newtheorem{example}[equation]{Example}

\theoremstyle{remark}

\newcommand{\beq}{\begin{equation}}
\newcommand{\eeq}{\end{equation}}

\DeclareMathOperator{\Hom}{{Hom}}

\DeclareMathOperator{\Tor}{Tor}

\DeclareMathOperator{\Ker}{Ker}
\DeclareMathOperator{\coker}{Coker}

\DeclareMathOperator{\gr}{gr}

\DeclareMathOperator{\Der}{Der}

\newcommand{\mc}{\mathcal}
\newcommand{\mf}{\mathfrak}
\newcommand{\mb}{\mathbb}

\renewcommand{\AA}{\mb A}

\newcommand{\CC}{\mb C}

\newcommand{\FF}{\mb F}
\newcommand{\GG}{\mb G}

\newcommand{\II}{\mb I}
\newcommand{\JJ}{\mb J}

\newcommand{\NN}{\mb N}

\newcommand{\WW}{\mb W}

\newcommand{\ZZ}{\mb Z}

 \newcommand{\kk}{{\Bbbk}}
 \newcommand{\Ua}{\mathrm{U}}

\DeclareMathOperator{\Rees}{Rees}

\newcommand{\normal}{\vartriangleleft}
\newcommand{\del}{\partial}
\newcommand{\wt}{\widetilde}
\newcommand{\wh}{\widehat}

\newcommand{\wtB}{\wt{B}}

\newcounter{onesidedW}[section]
\newcommand{\Wone}{
\ifthenelse{\value{onesidedW}=0}
    {\WW_1}
    {W_{\geq -1}}
}

\makeatletter
\newcommand{\subjclass}[2][1991]{%
  \let\@oldtitle\@title%
  \gdef\@title{\@oldtitle\footnotetext{#1 \emph{Mathematics Subject Classification.} #2}}%
}
\newcommand{\keywords}[1]{%
  \let\@@oldtitle\@title%
  \gdef\@title{\@@oldtitle\footnotetext{\emph{Key words and phrases.} #1.}}%
}
\makeatother

\title{On coherence of enveloping algebras of infinite-dimensional Lie algebras}
\author{Susan J. Sierra and James Timmins}
\date{\vspace{-3ex}}


\begin{document}
\renewcommand{\onlyinsubfile}[1]{} 
\renewcommand{\notinsubfile}[1]{#1}

\maketitle

\begin{abstract} \noindent
We investigate when the universal enveloping algebra $\Ua(\mf{g})$ of an infinite-dimensional Lie algebra $\mf{g}$ is coherent. We show that if $\mf{g}$ is the loop algebra of a finite-dimensional Lie algebra $\mf{h}$ then $\Ua(\mf{g})$ is coherent exactly when $\mf{h}$ is nilpotent. As a consequence, we show that $\Ua(\mf{g})$ is not coherent if $\mf{g}$ is complex, simple, $\mb{Z}$-graded, of polynomially bounded growth, and not isomorphic to the Witt algebra $W = \Der \CC[t,t^{-1}]$ or its one-sided version.
We also demonstrate that the homomorphic image of $\Ua(W)$ used to prove non-Noetherianity of $\Ua(W)$ in previous work of the first author and Walton is not coherent. On the other hand, we establish that two closely related connected graded rings are graded coherent. Finally, we prove that any exhaustively filtered ring is coherent if its Rees ring is graded coherent.
\end{abstract}

\tableofcontents

\subfile{Introduction.tex}

\subfile{Coherence_of_enveloping_algebras.tex}

\subfile{Graded_coherence_of_B_and_T.tex}

\subfile{Ungraded_coherence_and_Rees_rings.tex}

\bibliographystyle{amsalpha}
\bibliography{main.bib}

\end{document}

%% file: Introduction.tex
\section{Introduction}

A long standing open question on infinite-dimensional Lie algebras is: when is the universal enveloping algebra of a Lie algebra Noetherian? First appearing in the literature over half a century ago in \cite{amayo1974infinite}, and reformulated as a conjecture in \cite{sierra2014universal}, the answer to this question is believed to be: when the Lie algebra in question is finite-dimensional. Of course, any finite-dimensional Lie algebra has a Noetherian universal enveloping algebra as a natural consequence of the Poincaré–Birkhoff–Witt theorem. The converse statement remains open, and in full generality still  very difficult, although advances have been made in various classes and examples of infinite-dimensional Lie algebras, see \cite{sierra2014universal,andruskiewitsch2025noetherian,bell2025enveloping}.

In this article, we ask:

\begin{question} \label{Coherence of enveloping algebras question}
    When is the universal enveloping algebra of a Lie algebra \emph{coherent}?
\end{question}

A ring is left coherent when every finitely-generated left ideal is finitely-presented; this generalises being Noetherian. 
As universal enveloping algebras have a canonical anti-automorphism, for these rings the property of being coherent is left-right symmetric.
We also consider the property of being graded coherent, the obvious extension of the definition to graded ideals in graded rings. A graded ring which is graded coherent is not necessarily coherent, see \Cref{graded coherent not coherent example}. 

Every left Noetherian ring is automatically left coherent (since all finitely-generated modules are finitely-presented), though the class of coherent non-Noetherian rings contains many interesting examples. For example, commutative polynomials in infinitely many variables over any Noetherian ring, and any free algebra over a field (in more than one variable), are both coherent but non-Noetherian.

Coherent rings also extend the notion of Noetherianity in a representation-theoretic fashion: when $R$ is any ring, $R$ is left Noetherian if and only if the category of finitely-generated left $R$-modules is abelian. 
Analogously, $R$ is left coherent if and only if the category of finitely-\emph{presented} left $R$-modules is abelian.
An example where this is applied is the coherence of certain group rings \cite{shotton2020category,timmins2023coherence} implying categories of representations in the $p$-adic Langlands programme are \emph{locally coherent}, \cite[Conjecture 6.1.4]{emerton2025categorical}. For similar reasons, coherent rings also play an important role in noncommutative algebraic geometry, where one typically works with (a quotient category) of the graded module category over a graded ring, rather than the ring itself.
See \cite{polishchuk2005noncommutative,piontkovski2008coherent,smith}. 

We mention that the coherence of a ring can also be characterised homologically: a ring is left coherent if arbitrary products of flat right $R$-modules are flat, see \cite[Theorem 2.1]{chase1960direct}. This compares with the Bass-Papp theorem that a ring is left Noetherian if and only if any direct sum of injective left modules is injective.

Returning to \Cref{Coherence of enveloping algebras question}, observe that some potential answers can be quickly ruled out. If $\mf{g}$ is a finite-dimensional Lie algebra, then of course $\Ua(\mf{g})$ is coherent because it is Noetherian. The converse does not hold though, as if $\mf{a}$ is an infinite-dimensional abelian Lie algebra, $\Ua(\mf{a})$ is a polynomial ring with infinitely many variables, and is therefore coherent but non-Noetherian. On the other hand, not every Lie algebra has a coherent enveloping algebra -- in \Cref{shift non-coherent prop} we demonstrate that the shift algebra $\mf{s}$, an extension of a one-dimensional Lie algebra by a countable-dimensional abelian $\mf{a}$, has a non-coherent enveloping algebra.

Analogously to the conjectured connection between finite-dimensionality and Noetherianity, one might guess that coherence of universal enveloping algebras is controlled by the rate of growth in the Lie algebra. This possibility appears to be ruled out by the fact that the shift algebra $\mf{s}$ is finitely-generated and only of linear growth, but $\Ua(\mf{s})$ is not coherent. Conversely, if $F$ is any finitely-generated free Lie algebra, then $\Ua(F)$ is coherent, but $F$ has exponential growth (assuming it has more than one generator).

In this article, we consider the ``smallest'' infinite-dimensional Lie algebras, those of polynomial growth, and answer \Cref{Coherence of enveloping algebras question} for almost all cases when the Lie algebra is simple and $\mb{Z}$-graded.

\begin{restatable}{thm}{gradedLienoncoherence} \label{graded Lie algebras non-coherence thm}
    Let $\kk$ be an algebraically closed field of characteristic 0, and let 
    $\mf{g}$ be an infinite-dimensional $\mb{Z}$-graded simple $\kk$-Lie algebra with polynomially bounded growth which is not isomorphic to the Witt algebra $W = \Der \kk[t,t^{-1}]$ or the one-sided Witt algebra $\mb W_1 = \Der \kk[t]$. Then $\Ua(\mf{g})$ is not (graded) coherent. 
\end{restatable}

\Cref{graded Lie algebras non-coherence thm} relies on Mathieu's famous classification \cite{mathieu1992classification} of Lie algebras satisfying the hypotheses of the theorem, and is thus a case-by-case argument.
One case is when $\mf{g}$ is a loop algebra of a finite-dimensional simple Lie algebra. We deduce this case from a more general characterisation of coherence for loop algebras of finite-dimensional Lie algebras.

\begin{restatable}{thm}{loopcoherence} \label{loop coherence thm}
Let $\mf{g}$ be a finite-dimensional Lie algebra over any field $\kk$. The universal enveloping algebra $\Ua(L(\mf{g}))$ is a (graded) coherent ring if and only if $\mf{g}$ is nilpotent.
\end{restatable}

We further show that the enveloping algebra of any nilpotent Lie algebra is coherent, \Cref{nilpotent Lie algebra coherent prop}, which in turn implies one direction of \Cref{loop coherence thm}.

The remainder of the cases in \Cref{graded Lie algebras non-coherence thm}, namely the algebras of Cartan type, are proved to have non-coherent enveloping algebras by demonstrating that the shift algebra $\mf{s}$ is found as a (graded) subalgebra, see \Cref{Cartan shift algebra thm}.

The missing cases in \Cref{graded Lie algebras non-coherence thm} are the  Witt algebra $W$ and the  one-sided Witt algebra $\mb W_1$. 
These are also the difficult cases in the result \cite[Corollary~0.6]{sierra2014universal} that if $\mf{g}$ is an infinite-dimensional $\mb{Z}$-graded simple Lie algebra with polynomially bounded growth then $\Ua(\mf{g})$ is not Noetherian.  Thus 
\Cref{graded Lie algebras non-coherence thm} shows that, similarly to the question of Noetherianity, the coherence of universal enveloping algebras of $\mb{Z}$-graded Lie algebras of polynomial growth is determined by the answer for the Witt algebra.

Unfortunately, we have been unable to determine if the universal enveloping algebra of the Witt algebra or the one-sided Witt algebra is a coherent ring. However, during the course of our efforts to answer this question, we have determined (graded) coherence and non-coherence of some algebras closely related to these enveloping algebras.

To describe these, we recall two important aspects of the proof of non-Noetherianity of the Witt algebra $W$. One is the first \emph{orbit homomorphism} (as named  in \cite{phamXorbit})
\[ \Psi: \Ua(W) \rightarrow \widetilde{A_1} \otimes_{\kk} \kk[v_0], \]
defined in \cref{Non-coherence of B}, a generalisation of the natural inclusion of the Witt algebra into the localised Weyl algebra 
\[ \wt{A_1} = A_1[t^{-1}] = \kk[t,t^{-1}][\del; \frac{d}{dt}]. \] 
The essence of the proof that $\Ua(W)$ is non-Noetherian in \cite{sierra2014universal} is to show that the image of $\Psi$ is non-Noetherian, from which the result follows immediately.

In this article, we go further, and demonstrate that the image of $\Psi$ is not even coherent.
Note that $\Wone$ is isomorphic to a subalgebra of $W$ and so $\Psi$ restricts to a homomorphism from $U(\Wone) \to A_1 \otimes \kk[v_0]$.

\begin{restatable}{thm}{Bnotcoherent} \label{B not coherent thm}
    The images of $\Ua(W)$ and $\Ua(\Wone)$ under $\Psi$ are not left or right (graded) coherent rings.
\end{restatable} 

It is still possible that one or both of the enveloping algebras $\Ua(W)$, $\Ua(\WW_{1})$ may be coherent, since not every quotient of a coherent ring is coherent. (Note that if $\Ua(W)$ is coherent, so is $\Ua(\WW_{1})$, see \Cref{coherence descent lemma}.)
It is known that if $R$ is a coherent ring and $I$ is a two-sided ideal finitely-generated as a left ideal, then $R/I$ is coherent, see \Cref{finitely-presented coherence prop}. However, this does not apply in our situation, because the kernel of $\Psi$ is not finitely-generated as a left ideal by \cite{sierra2016maps}, see also \cite[Theorem 1.4]{phamXkernel}.

We refine our results further, by considering a second important aspect of \cite{sierra2014universal}. The proof of non-Noetherianity of $\Psi(\Ua(W))$, and hence $\Ua(W)$, was inspired by observing that the image of the orbit homomorphism $\Psi$ lies in a geometrically motivated idealiser, called $\wt{T}$ here. The ring $\wt{T}$ has the interesting property that it is left Noetherian but not right Noetherian.

In more detail, $\wt{T}$ is a graded subring of
\[\wt{S} =\wt{A_1} \otimes_{\kk} \kk[v_0] \cong \kk[a,b][t,t^{-1}; \sigma] \]
that contains $\wtB = \Psi(\Ua(W))$. More precisely, $\wt{T}$ is the idealiser of the left ideal $\wt{S}a+ \wt{S}b$ of $\wt{S}$, and is given by $\wt{T} = \kk \oplus (
\wt{S}a+\wt{S}b)$, see \Cref{T def}.
Extending the study of $\wt{T}$, and contrasting with \Cref{B not coherent thm}, we prove that $\wt{T}$ and its subring
\beq\label{T} T = \wt{T} \cap \kk[a,b][t;\sigma] \eeq
are graded right coherent.

\begin{restatable}{thm}{Tgradedcoherent} \label{T graded coherent thm}
    Although $T$ and $\wt{T}$ are not right Noetherian, both are left Noetherian and graded right coherent.
\end{restatable}

That $\wt{T}$ is graded coherent and $\wtB$ is not demonstrates the subtlety of coherence, and indicates that further work is necessary to offer a hypothesis on the coherence of $\Ua(W)$.

Since we are unable to resolve the question of coherence of $\Ua(W)$ any further, we turn our attention to the \emph{un}graded coherence of $T$ and $\wt{T}$. We prove a useful result for resolving questions of coherence via graded coherence, as follows.

\begin{restatable}{thm}{Reescoherence} \label{Rees coherence thm}
    Let $A$ be a filtered ring. If the Rees ring of $A$ is graded left coherent, then $A$ is left coherent.
\end{restatable}

Recall that $A$ (under a complete filtration) is Noetherian if its associated graded ring is (graded) Noetherian. This is false for coherence: any non-coherent enveloping algebra (for example those of \Cref{graded Lie algebras non-coherence thm}) with its PBW filtration provides a counterexample. \Cref{Rees coherence thm} partially rectifies this failure.

Now, $T$ is a graded ring, and therefore may be made a filtered ring by taking increasing unions of graded pieces. Suspecting that $T$ is coherent, we consider the Rees ring
\[ \mc{R} = \Rees(T) \subseteq T[s],\]
and the Rees ring $\wt{\mc{R}}$ of $\wt{T}$, and study their graded coherence. Unfortunately, we are unable to prove that $\wt{\mc{R}}$ or $\mc{R}$ is graded right coherent under the canonical grading coming from the filtration on $T$, realised by the $s$-degree.

However, the original grading on $\wt{T}$ also induces a grading on $\wt{\mc{R}}$ by $ts$-degree, which we call the $u$-grading. This grading gives $\wt{\mc{R}}$ a geometric flavour similar to $\wt{T}$, and we obtain graded right coherence of $\mc{R}$ and $\wt{\mc{R}}$ under this grading.

\begin{restatable}{thm}{ugradedcoherence} \label{u-graded coherence thm}
    The rings $\mc{R}$ and $\wt{\mc{R}}$ are graded right coherent with respect to the $u$-grading.
\end{restatable}

In fact, $\wt{\mc{R}}$ is bi-graded by these two gradings, and an immediate consequence of \Cref{u-graded coherence thm} is that $\mc{R}$ (and $\wt{\mc{R}}$) is bi-graded right coherent, meaning that every finitely-generated bi-graded right ideal is a finitely-presented $\mc{R}$-module, \Cref{bi-graded coherence cor}.

Our results touch on but do not answer a number of open questions surrounding coherence of enveloping algebras of the Witt algebras $W$, $\Wone$. The most immediate are as follows.

\begin{question}
    \begin{enumerate}[label=(\alph*)]
        \item Are $\Ua(W)$, $\Ua(\WW_{1})$ coherent rings?
        \item Are $\Ua(W)$, $\Ua(\WW_{1})$ graded coherent?
    \end{enumerate}
\end{question}

Whilst we would guess that parts (a) and (b) of the above question have the same answer, we make no conjecture as to the answer, due to the contradictory indications we receive from \Cref{B not coherent thm} and \Cref{T graded coherent thm}. It is also possible that $\Ua(\WW_{1})$ could be (graded) coherent whilst $\Ua(W)$ is not.

\begin{question}
    Are the graded right coherent rings $T$ and $\wt{T}$ right coherent?
\end{question}

If $T$ or $\wt{T}$ is not right coherent, it would be the first ``natural" example to our knowledge of a graded coherent ring which is not coherent. 
Of course, despite \Cref{Rees coherence thm} it is possible that $\mc{R} = \Rees(T)$ is not graded coherent whilst $T$ is coherent.

\begin{question}
    Are the Rees rings $\mc{R}$, $\wt{\mc{R}}$ graded right coherent or not?
\end{question}

We briefly discuss the structure of the paper.  
\Cref{graded Lie algebras non-coherence thm} is proved in \cref{Coherence of graded simple Lie algebras of polynomial growth}, with \Cref{loop coherence thm} preceding it in \cref{Coherence of loop algebras}. In \Cref{Graded noncoherence of some rings related to UW} we study rings derived from the enveloping algebra of the Witt algebra. \Cref{B not coherent thm} is proved in \cref{Non-coherence of B}, and we prove \Cref{T graded coherent thm} in \cref{Graded coherence of T subsection}. In \Cref{Ungraded coherence and Rees rings section}, we turn to questions of un-graded coherence, and after demonstrating \Cref{Rees coherence thm} in \cref{Coherence of Rees rings subsection}, we then prove \Cref{u-graded coherence thm} in \cref{Proof of u-graded coherence thm subsection}.

\emph{Acknowledgements.} Both authors were supported by EPSRC grant EP/T018844/1. The second author was additionally supported by the ERC under the EU Horizon 2020 research and innovation programme, under grant agreement 948885, and by the Royal Society URF/R/261036.

\emph{Generative AI use.} No generative AI tools were used in the preparation of this manuscript.

%% file: Coherence_of_enveloping_algebras.tex
\section{Background and notation}

Throughout we fix a ground field $\kk$.   
All vector spaces, (Lie) algebras, linear maps, etc. will be defined over $\kk$ and unadorned $\otimes$ will mean $\otimes_{\kk}$.
We denote by $\NN  = \{ n \in \ZZ \mid n \geq 0\}$ the natural numbers.

If $A = \bigoplus_{n \in \NN} A_n$ is an $\NN$-graded algebra,
we say $A$ is {\em finitely graded} or {\em finitely $\NN$-graded} if $\dim A_n < \infty$ for all $n \in \NN$, and {\em connected graded} if it is finitely graded and $A_0 = \kk$.

A ring is {\em left coherent} if every finitely-generated left ideal is finitely-presented. Similarly, a graded ring is {\em graded left coherent} if every finitely-generated graded left ideal is finitely-presented.
The right-hand notions have the obvious definitions, mutatis mutandis.
We will often state symmetric results only for one side for conciseness. We say a ring is (graded) coherent when it is both left and right (graded) coherent.

We now give some fundamental results on coherent rings.

The following is \cite[Theorem 2.3.3]{glaz1989commutative}, with the proof being identical for commutative or non-commutative rings, and the graded extension being straightforward.

\begin{lemma} \label{direct limit coherence lemma}
    Let $R$ be a ring which is a directed colimit of left (graded) coherent subrings $S$ (along graded ring homomorphisms), such that $R$ is a flat right $S$-module for each $S$. Then $R$ is left (graded) coherent.
\end{lemma}

The following two results can be found in \cite{harris1966some1} for ungraded rings, with the proof for graded versions being identical.

\begin{proposition}[{\cite[Corollary 2.1]{harris1966some1}}] \label{faithfully flat coherence prop}
    Let $\phi: A \rightarrow B$ be a (graded) ring homomorphism. If $B$ is (graded) left coherent and faithfully flat as a right $A$-module, then $A$ is (graded) left coherent.
\end{proposition}

\begin{proposition}[{\cite[Corollary 1.2]{harris1966some1}}] \label{finitely-presented coherence prop}
    Let $\phi: A \rightarrow B$ be a (graded) ring homomorphism. If $A$ is (graded) left coherent and $B$ is finitely-presented as a left $A$-module, then $B$ is (graded) left coherent. 
\end{proposition}

As a consequence, coherence of a $\kk$-algebra can be reduced to the case where $\kk$ is algebraically closed.

\begin{lemma} \label{algebraic closure coherent lemma}
    Let $A$ be a $\kk$-algebra and let $\bar{\kk}$ be the algebraic closure of $\kk$.
    Then $A$ is left coherent if and only if $\bar{\kk} \otimes A$ is left coherent.
\end{lemma}
\begin{proof}
    Suppose $A$ is left coherent. Let $\ell$ be a finite field extension of $\kk$. Then $\ell \otimes A = A \otimes \ell$ is a finitely-generated free left $A$-module, hence is left coherent by \Cref{finitely-presented coherence prop}. Then the direct limit over all finite extensions
    \[ \varinjlim_{\ell} (\ell \otimes A) = \bar{\kk} \otimes A\]
    is left coherent by \Cref{direct limit coherence lemma}.
    
    On the other hand, $\bar{\kk} \otimes A$ is a free left $A$-module, so the converse follows from \Cref{faithfully flat coherence prop} via faithful flatness.
\end{proof}

If $A$ is a graded algebra with $\kk$ in degree zero, then the result of \Cref{algebraic closure coherent lemma} also holds for graded coherence.

\begin{proposition} \label{fg intersection prop}
    Let $R$ be a domain. Then $R$ is left coherent if and only if the intersection of any two finitely-generated left ideals is finitely-generated.
\end{proposition}
\begin{proof}
Since $R$ is a domain, this follows directly from Theorem 2.3.2(7) of \cite{glaz1989commutative}, noting that the result is stated for commutative rings but the proof holds in the noncommutative case.    
\end{proof}

It is well known that an $\NN$-graded ring which is graded Noetherian must be Noetherian. However, this does not extend to (graded) coherence.

\begin{example} \label{graded coherent not coherent example}
    Let $R$ be a left coherent ring such that the polynomial ring $R[X]$ is not left coherent. 
    Such rings exist: an example due to Soublin is given in \cite[7.3.13]{glaz1989commutative}. Give $R[X]$ the grading by $X$-degree, where $X$ is of degree 1 and all elements of $R$ are of degree zero. Then by \cite[Theorem 1.1]{minamoto2018criterion}, $R[X]$ is graded left coherent but not left coherent. 
\end{example}

\Cref{graded coherent not coherent example} gives the only examples known to us of rings which are graded coherent but not coherent. 
In \cref{Coherence results for the Rees ring of T} we investigate the coherence of the graded coherent rings $\wt{T}$ and $T$, defined in \eqref{T}.

\section{Coherence of enveloping algebras}

If $\mf g$ is a finite-dimensional Lie algebra, then $\Ua(\mf g)$ is Noetherian and so necessarily coherent.
The question of coherence of $\Ua(\mf g)$  is much more subtle and interesting if $\mf g$ is infinite-dimensional.  
For example, if $\mf g$ is abelian then $\Ua(\mf g)$ is a direct limit of finitely generated polynomial rings and so is coherent; if $\mf g$ is free then $\Ua(\mf g)$ is also coherent.

However, we find in this section that enveloping algebras of many infinite-dimensional Lie algebras are not coherent: for example, loop algebras of simple Lie algebras. In fact we prove that if $\mf g$ is a finite dimensional Lie algebra and $L(\mf g)$ is its loop algebra, then $\Ua(L(\mf g))$ is coherent if and only if $\mf g$ is nilpotent, \Cref{loop coherence thm}.
We thereby prove that if $\mf g$ is any infinite-dimensional graded simple Lie algebra of polynomial growth which is not isomorphic to a Witt algebra, then $\Ua(\mf g)$ is not coherent, \Cref{graded Lie algebras non-coherence thm}. Unfortunately we are not able to resolve the coherence question for $\Ua(\WW_1)$ or $\Ua(W)$.

Our main technique to prove \Cref{graded Lie algebras non-coherence thm} and \Cref{loop coherence thm} is the observation that since $\Ua(\mf s)$ is not coherent, see \Cref{shift non-coherent prop} where $\mf{s}$ is the {\em shift algebra} defined below, a sufficient condition for $\Ua(\mf g)$ to not be coherent is that $\mf s$ is a subalgebra of $\mf g$. However, $\mf s$ is not a subalgebra of $\WW_1$ or of $W$, and so we are able to say less about those enveloping algebras. We investigate the coherence of some rings related to $\Ua(\WW_1)$ and $\Ua(W)$ in the next section.

In this section, the word ``graded" may be dropped from all results mentioning ``graded coherent", with the proofs remaining identical.

\subsection{The shift algebra}

\begin{defn} \label{shift algebra def}
    The {\em shift algebra} $\mf{s}$ is the infinite-dimensional Lie algebra  with basis $\{h, e_1, e_2, \dots\}$ and Lie brackets given by $[h, e_i]=e_{i+1}$, $[e_i, e_j]=0$ for all $i,j \geq 1$.
\end{defn}

We give the shift algebra the grading where $\deg(h)=1$ and $\deg(e_i)=i$.

We give a simple lemma which allows us to more easily recognise the shift algebra.

\begin{lemma} \label{recognise shift lemma}
    Let $\mf{g}$ be a Lie algebra with a $\kk$-basis $\{H, E_1, E_2, \dots\}$ such that $[H, E_i]\in \kk^\times E_{i+1}$ and $[E_i, E_j]=0$ for all $i, j \geq 1$. Then $\mf{g} \cong \mf{s}$.
\end{lemma}
\begin{proof}
    Let $[H, E_i] = \lambda_i E_{i+1}$ for all $i \geq 1$. Then the $\kk$-linear map $\phi: \mf{s} \rightarrow \mf{g}$ with $\phi(h)=H$ and
    \[ \phi(e_i) = \lambda_1\dots\lambda_{i-1}E_i,\]
    is an isomorphism of $\kk$-Lie algebras.
\end{proof}

\begin{proposition} \label{shift non-coherent prop}
    The universal enveloping algebra $\Ua(\mf{s})$ is not a (graded) coherent ring.
\end{proposition}

Although \Cref{shift non-coherent prop} can be proved directly, we give a more general result based on the presentation of an (associative) $\kk$-algebra.

\begin{proposition} \label{presentation of algebra non-coherence prop}
    Let $A$ be a connected $\NN$-graded $\kk$-algebra, finitely-generated by homogeneous elements of positive degree. If $A$ is not finitely-presented as a $\kk$-algebra, then $A$ is not (graded) coherent. 
\end{proposition}
\begin{proof}
    Because $A$ is finitely-generated by homogeneous elements, $A$ is isomorphic to a quotient of a weighted free algebra $F$ by a homogeneous two-sided ideal $I$. Let $J$ be the ideal of $A$ generated by all positive degree homogeneous elements, so $A/J \cong \kk$. By Lemma 2.1.3 of \cite{rogalski}, $J$ is finitely-generated as a left ideal, and if $J$ is finitely-presented, then $I$ is a finitely-generated two-sided ideal of $F$. Hence, if $A$ is not a finitely-presented algebra, then $I$ cannot be finitely-generated, hence $A$ is not graded left coherent. The same reasoning holds for right coherence.
\end{proof}

\begin{lemma} \label{fp Lie algebra lemma}
    Let $\mf{g}$ be a  Lie algebra.  
    Then $\mf{g}$ is finitely-presented as a Lie algebra if and only if $\Ua(\mf{g})$ is a finitely-presented (associative) algebra.
\end{lemma}
\begin{proof}
    Let $F$ be a free Lie algebra and $R$ be a Lie ideal of $F$. Let $i: F \rightarrow \Ua(F)$ be the natural inclusion map, and let  $I_R = \Ua(F)i(R)\Ua(F)$. The Lie algebra homomorphism $F/R \rightarrow \Ua(F)/I_R$ extends to a $\kk$-algebra homomorphism
    $\Ua(F/R) \rightarrow \Ua(F)/I_R$, by the universal property of universal enveloping algebras. On the other hand, the map of $\kk$-algebras $\Ua(F) \rightarrow \Ua(F/R)$ induced by $F \rightarrow F/R$ has kernel containing $I_R$, from which we obtain an isomorphism
    \[  \Ua(F)/I_R\cong \Ua(F/R).\]
    For any subset $S \subseteq F$ generating the Lie ideal $\langle S \rangle$ of $F$, we have $I_{\langle S \rangle} = \Ua(F) i(S) \Ua(F)$, due to the inclusions
    \[ i(\langle S \rangle) \subseteq \Ua(F) i(S) \Ua(F) \subseteq I_{\langle S \rangle}. \]
    Hence, $I_R$ is a finitely-generated two-sided ideal if and only if there is a finite subset $S$ with $\langle S \rangle = R$, if and only if $R$ is finitely-generated as a Lie ideal.
    
    If $\mf{g}$ is finitely-presented, we have $\mf{g} \cong F/R$ where $F$ is a finitely-generated free Lie algebra, and $R$ is a finitely-generated Lie ideal. Then $\Ua(F)$ is a finitely-generated free $\kk$-algebra and $I_R$ is a finitely-generated two-sided ideal. Hence $\Ua(\mf{g})$ is finitely-presented.
    
    Conversely, if $\Ua(\mf{g})$ is finitely-presented, then it is generated by finitely many elements of $\mf{g}$.
    This gives a surjection $F \rightarrow \mf{g}$ of Lie algebras where $F$ is a finitely-generated free Lie algebra. 
    Let $R$ be the kernel of this homomorphism so that $\mf{g} \cong F/R$ and $\Ua(\mf{g}) \cong \Ua(F)/I_R$.
    Then $\Ua(F)$ is a finitely-generated free $\kk$-algebra, so $I_R$ must be finitely-generated since $\Ua(\mf{g})$ is finitely-presented, implying that $R$ is finitely-generated as a Lie ideal. Thus $\mf{g}$ is finitely-presented.
\end{proof}

\begin{proof}[Proof of \Cref{shift non-coherent prop}]
    The shift algebra $\mf{s}$ is generated by two elements $h$ and $e_1$. Notice that the shift algebra has an alternative grading to that in \Cref{shift algebra def}, where $h$ has degree $1$, and $e_j$ degree $j+1$. 
    In this grading, $\mf{s}$ is {\em filiform}; such Lie algebras are studied and classified in \cite{fialowski1983classification}.
   By \cite[Corollary~1]{fialowski1983classification}  $\mf{s}$ is not finitely-presented. 
   (We note that \cite{fialowski1983classification} has an assumption that the ground field has characteristic $0$, but the argument that $\mf{s}$ is not finitely-presented is characteristic-free.)
   By \Cref{fp Lie algebra lemma}, $\Ua(\mf{s})$ is not finitely-presented as a $\kk$-algebra. 
   It then follows from \Cref{presentation of algebra non-coherence prop} that $\Ua(\mf{s})$ is not graded coherent, under the usual grading.
\end{proof}

\subsection{Coherence of loop algebras} \label{Coherence of loop algebras}

\begin{defn} \label{loop algebra def}
    Let $\mf{g}$ be a Lie algebra. The {\em loop algebra} $L(\mf{g})$ is the vector space $\mf{g}[t,t^{-1}]$, with Lie bracket given by the unique $\kk$-linear extension of $[xt^m, yt^n] = [x,y]t^{m+n}$ for all $x,y \in \mf{g}$, $m,n \in \mb{Z}$.
\end{defn}

The loop algebra of $\mf{g}$ has a natural $\ZZ$-grading with the $n$th graded piece being $\mf{g}t^n$.
In this section, we consider the coherence of  enveloping algebras of loop algebras.

\begin{lemma} \label{coherence descent lemma}
    Let $\mf{g}$ be a (graded) Lie algebra and $\mf{h}$ be a (graded) subalgebra. If $\Ua(\mf{g})$ is (graded) coherent, then $\Ua(\mf{h})$ is (graded) coherent.
\end{lemma}
\begin{proof}
    The enveloping algebra $\Ua(\mf{g})$ is a free, hence faithfully flat, right $\Ua(\mf{h})$-module. The result follows from \Cref{faithfully flat coherence prop}.
\end{proof}

\begin{defn} 
Let $\mf{b}$ be the two-dimensional solvable Lie algebra, that is $\mathfrak b$ has  basis $\{x,y\}$ such that $[x,y] = y$.
\end{defn}

\begin{corollary} \label{two dim subalgebra non-coherence cor}
Let $\mf{g}$ be a Lie algebra with a subalgebra isomorphic to $\mf{b}$. Then $\Ua(L(\mf{g}))$ is not graded coherent.
\end{corollary}
\begin{proof}
    Notice that $\mf{s}$ is isomorphic to a graded subalgebra of the loop algebra $L(\mf{b})$, by identifying $h$ with $xt$ and $e_n$ with $yt^n$. Since $\Ua(\mf{s})$ is not (graded) coherent by \Cref{shift non-coherent prop}, it follows from \Cref{coherence descent lemma} that $\Ua(L(\mf{b}))$ is not (graded) coherent. Clearly, $L(\mf{b})$ is a graded subalgebra of $L(\mf{g})$, so another application of \Cref{coherence descent lemma} gives the result.
\end{proof}

\begin{proposition} \label{nilpotent Lie algebra coherent prop}
    Let $\mf{g}$ be a nilpotent Lie algebra. Then $\Ua(\mf{g})$ is a coherent ring.
\end{proposition}
\begin{proof}
    If $\mf{g}$ is a finitely-generated nilpotent Lie algebra, then $\mf{g}$ is finite-dimensional because only finitely many iterated Lie brackets in the generators are non-zero. Thus, $\Ua(\mf{g})$ is Noetherian. Then if $\mf{g}$ is an arbitrary nilpotent Lie algebra, $\mf{g}$ is a direct limit of its finitely-generated subalgebras, hence is a direct limit of finite-dimensional subalgebras. Therefore $\Ua(\mf{g})$ is a direct limit of Noetherian subalgebras, and is flat over each subalgebra. Thus $\Ua(\mf{g})$ is coherent by \Cref{direct limit coherence lemma}. 
\end{proof}

\begin{corollary} \label{nilpotent loop coherence prop}
Suppose $\mf{g}$ is nilpotent. Then $\Ua(L(\mf{g}))$ is a coherent ring.
\end{corollary}
\begin{proof} 
Due to \Cref{nilpotent Lie algebra coherent prop}, it is enough to show that $L(\mf{g})$ is nilpotent. But, it is straightforward to check that for any subalgebra $\mf{h}$ of $\mf{g}$,
\[ [L(\mf{g}), L(\mf{h})] = L([\mf{g}, \mf{h}]),\]
from which the result follows.
\end{proof}

Using \Cref{two dim subalgebra non-coherence cor}, we will   prove the converse of \Cref{nilpotent loop coherence prop}.

\begin{lemma} \label{non-nilpotent subalgebra lemma}
Suppose $\mf{g}$  is a finite-dimensional Lie algebra that is not nilpotent and $\kk$ is algebraically closed. Then $\mf{g}$ has a Lie subalgebra isomorphic to $\mf{b}$.
\end{lemma}
\begin{proof} 
Since $\mf{g}$ is finite-dimensional, Engel's theorem implies that there exists $x \in \mf{g}$ such that 
\[ \mathrm{ad}(x): \mf{g} \rightarrow \mf{g} \]
is a non-nilpotent linear map. Since $\kk$ is algebraically closed, $\mathrm{ad}(x)$ has at least one eigenvector and eigenvalue. Not all eigenvalues of $\mathrm{ad}(x)$ are zero, otherwise by consideration of Jordan normal form, $\mathrm{ad}(x)$ is nilpotent, a contradiction. So there exists $\lambda \in \kk^\times$ and $y \in \mf{g}\backslash\{0\}$ such that
\[ [x,y] = \lambda y. \]
By considering the elements $\lambda^{-1}x, y \in \mf{g}$, it follows that the subalgebra of $\mf{g}$ with basis $\{x,y\}$ is isomorphic to $\mf{b}$. 
\end{proof}

\notinsubfile{ 
\loopcoherence*
}

\onlyinsubfile{
\begin{thm} \label{loop coherence cor}
Let $\mf{g}$ be a finite-dimensional Lie algebra. The universal enveloping algebra $\Ua(L(\mf{g}))$ is a (graded) coherent ring if and only if $\mf{g}$ is nilpotent.
\end{thm}
}

\begin{proof}
Since $\bar{\kk} \otimes \mf{g}$ is nilpotent if and only if $\mf{g}$ is nilpotent, and $\Ua(L(\bar{\kk} \otimes \mf{g})) = \bar{\kk} \otimes \Ua(L(\mf{g}))$, we may without loss of generality assume $\kk$ is algebraically closed, by \Cref{algebraic closure coherent lemma}.

The case that $\mf{g}$ is nilpotent is shown by \Cref{nilpotent loop coherence prop}. If $\mf{g}$ is not nilpotent, then by \Cref{non-nilpotent subalgebra lemma}, $\mf{g}$ has a subalgebra isomorphic to $\mf{b}$, and thus by \Cref{two dim subalgebra non-coherence cor}, $\Ua(L(\mf{g}))$ is not (graded) coherent. 
\end{proof}

Finally, we consider a related class of Lie algebras: the \emph{twisted} loop algebras. For the remainder of this subsection, specialise $\mf{g}$ to be finite-dimensional and simple and $\kk$ to be algebraically closed of characteristic zero. Let $\sigma$ be a Lie algebra automorphism of $\mf{g}$ of finite order $r$ and assume $\kk$ contains a primitive $r$th root of unity $\omega$. The action of $\sigma$ extends to $L(\mf{g})$ via $\sigma(gt^n) = \sigma(g)\omega^{-n}t^n$, and the twisted loop algebra is $L(\mf{g})^\sigma$, the $\sigma$-invariants in $L(\mf{g})$.

\begin{thm} \label{twisted loop algebra non-coherence thm}
    The enveloping algebra of the twisted loop algebra $L(\mf{g})^\sigma$ is not graded coherent.
\end{thm}
\begin{proof}
As above, we demonstrate that $L(\mf{g})^\sigma$ contains a copy of the shift algebra $\mf{s}$ as a subalgebra. We refer to subsection 2.1 of \cite{biswal2025ideals} for basic facts on twisted loop algebras.

Indeed, $\sigma$ corresponds, without loss of generality, to a Dynkin diagram automorphism of order $r$. Then \cite[2.1(0), p1203]{biswal2025ideals} shows that the subalgebra generated by homogeneous elements of degree divisible by $r$ is equal to $\mf{g}^\sigma[t^r, t^{-r}]$, which is isomorphic to $L(\mf{g}^\sigma)$ up to a rescaling of the grading. Thus by \Cref{coherence descent lemma} and \Cref{loop coherence thm} it is enough to check that $\mf{g}^\sigma$ is not nilpotent. But this follows from \cite[2.1(7), p1203-1204]{biswal2025ideals}.
\end{proof}

\begin{corollary} \label{affine LA non-coherence}
The universal enveloping algebra of an affine or affine Kac-Moody Lie algebra is not coherent.
\end{corollary}
\begin{proof}
Let $\wt{\mf g}$ be an affine Kac-Moody algebra.
This means that there is a finite-dimensional simple Lie algebra $\mf g$, 
a twisted loop algebra $L(\mf g)^\sigma$ of $\mf g$,
and a central extension $\widehat{\mf g}$ of $L(\mf g)^\sigma$ by a central element $c$ so that $[\wt{\mf g}, \wt{\mf{g}}] \cong \widehat{\mf g}$.
(We note that the affine Lie algebras are precisely the algebras of the form $\widehat{\mf g}$.)
Further, there is $d \in \wt{\mf g}$ so that $\wt{\mf g} \cong \widehat{\mf g} \rtimes \kk d$.
Thus the adjoint action of $d$  induces a derivation $\delta$ of  $\Ua(\widehat{\mf g}) $, 
and $\Ua(\wt{\mf g}) \cong \Ua(\widehat{\mf g})[d;\delta]$.

By \Cref{twisted loop algebra non-coherence thm} $\Ua(L(\mf g)^\sigma) \cong \Ua (\widehat{\mf g})/(c)$ is not coherent,
so by \Cref{finitely-presented coherence prop} $\Ua(\widehat{\mf g})$ is not coherent.
Thus by \Cref{faithfully flat coherence prop} $\Ua(\wt{\mf g})$ is not coherent.
\end{proof}

\subsection{Coherence of graded simple Lie algebras of polynomial growth} \label{Coherence of graded simple Lie algebras of polynomial growth}

In this subsection we apply the results of the previous two subsections to enveloping algebras of $\mb{Z}$-graded simple Lie algebras of polynomial growth. 
In this subsection, we suppose that $\kk$ is an algebraically closed field of characteristic zero.

Our proof is based on the classification of these Lie algebras, which was achieved by Mathieu in \cite{mathieu1992classification}. 

\begin{thm}
\label{graded Lie algebras classification prop}
    Every $\mb{Z}$-graded simple Lie algebra over $\kk$ with polynomially bounded growth is isomorphic to (exactly) one of the following.
    \begin{itemize}
        \item A finite-dimensional simple Lie algebra $\mf{h}$
        \item A (possibly twisted) loop algebra $L(\mf{h})$ with $\mf{h}$ as above
        \item A Lie algebra of Cartan type $\mb{W}_n$, $\mb{S}_n$, $\mb{H}_{2n}$, $\mb{K}_{2n+1}$
        \item The Witt algebra $W$.
    \end{itemize}
\end{thm} 

We now describe the Lie algebras in \Cref{graded Lie algebras classification prop}.
The {\em Witt algebra} is $ W = \Der(\kk[t,t^{-1}])$,
    where $\Der(A)$ denotes the $\kk$-linear derivations of a $\kk$-algebra $A$.
The {\em  Cartan algebras  } are $\mb W_n = \Der \kk[x_1, \dots, x_n]$, for $n \geq 1$, and the other algebras of Cartan type are certain subalgebras of the Cartan algebras, defined below.

Note that $\WW_1$ may be identified with a subalgebra of $W$.  More concretely,  if we let $\del = d/dt$, then $W$ has basis $\{ e_n = t^{n+1} \del \mid n \in \ZZ\}$, and is $\ZZ$-graded if we give $e_n$ degree $n$.
Then $\WW_1$ is isomorphic to $W_{\geq -1} := \kk\{ e_n \mid n \geq -1\}$, a subalgebra of $W$. We refer to $\WW_1 \cong W_{\geq-1} $ as the {\em one-sided Witt algebra}.

The main result of this section is:

\notinsubfile{
\gradedLienoncoherence*
}

\onlyinsubfile{
\begin{thm} \label{graded Lie algebras non-coherence thm}
    Let $\mf{g}$ be an infinite-dimensional $\mb{Z}$-graded simple Lie algebra with polynomially bounded growth, and which is not isomorphic to  $W$ or  $\mb W_1$. Then $\Ua(\mf{g})$ is not (graded) coherent. 
\end{thm}
}

We do not know if $\Ua(W)$ or $\Ua(\Wone)$ are (graded) coherent.
Note that if $\Ua(W)$ is (graded) coherent, then so is $\Ua(\Wone)$, by \Cref{coherence descent lemma}.

\begin{proposition} \label{Cartan shift algebra thm}
    Let $L$ be a Lie algebra of Cartan type with $L \neq \mb{W}_1$. Then $L$ contains a graded subalgebra isomorphic to the shift algebra $\mf{s}$.
\end{proposition}
\begin{proof}
    We use the description of the Cartan algebras found in \cite{dusautoy2004motivic}: let $n \geq 2$ and recall that $\mb{W}_n$ is the derivations of the polynomial ring $P_n = \kk[x_1, \dots, x_n]$.
    Let $D_i$ be the partial derivative with respect to $x_i$, and let $X = x_1^2 D_1$ and $Y_i = x_1^iD_2$ for all $i \geq 1$.
    Then it is straightforward to compute
    \[ [X, Y_i] = iY_{i+1}, \quad [Y_i, Y_j]=0,\]
    for all $i, j \geq 1$. Thus by \Cref{recognise shift lemma}, the Lie subalgebra of $\mb{W}_n$ generated by $\{X, Y_i \mid i \geq 1\}$ is graded and isomorphic to the shift algebra $\mf{s}$.
    
    The special subalgebra $\mb{S}_n$ of $\mb{W}_n$ is given by
    \[ \mb{S}_n = \Big\{ \sum_{j=1}^n p_jD_j \mid \sum_{j=1}^n D_j(p_j)=0 \Big\}.\]
    Again define $Y_i = x_1^iD_2$, and let $X = x_1^2D_1 - 2x_1x_2D_2 \in \mb{S}_n$. Then
    \[ [X, Y_i] = (i-2)Y_{i+1}, \quad [Y_i, Y_j]=0,\]
    for all $i, j \geq 1$, thus by \Cref{recognise shift lemma}, the subalgebra of $\mb{S}_n$ generated by $\{X, Y_i, \mid i \geq 3 \}$ is graded and isomorphic to the shift algebra.
    
    Now suppose that $n$ is even and let $n=2m$. The Hamiltonian subalgebra $\mb{H}_{n}$ consists of elements $D_H(p)$ for all $p \in P_n$, where
    \[ D_H(p) = \sum_{j=1}^{m} D_j(p)D_{m+j} - D_{m+j}(p)D_j.\]
    Then
    \[ Y_i = x_1^iD_{m+1}, \quad X = D_H(x_1^2x_{m+1})= -x_1^2D_{1} + 2x_1x_{m+1}D_{m+1},\]
    are elements of $\mb{H}_{n}$, and
    \[ [X, Y_i] = (2-i) Y_{i+1}, \quad [Y_i, Y_j]=0,\]
    so the subalgebra of $\mb{H}_n$ generated by $\{X, Y_i \mid i \geq 3\}$ is graded, and isomorphic to the shift algebra by \Cref{recognise shift lemma}.
    
    Finally, suppose $n$ is odd and let $n=2m+1$. The contact subalgebra $\mb{K}_n$ consists of elements $D_K(p)$ for all $p \in P_n$, where (expanding the description in \cite{dusautoy2004motivic})
    \begin{align*} D_K(p) &= D_{2m+1}(p)\sum_{j=1}^{2m}  x_jD_j + \sum_{j=1}^m D_j(p)D_{m+j}-D_{m+j}(p)D_j\\ &+ \Big( 2p - \sum_{j=1}^m x_jD_j(p) + x_{m+j}D_{m+j}(p)\Big)D_{2m+1}.\end{align*}
    The contact subalgebra has the grading induced from the polynomial ring $P_{2m+1}$ where $x_{2m+1}$ has degree 2, whilst $x_i$ has degree $1$ if $i \leq 2m$.
    We can compute that
    \[ D_K(x_1^i) = ix_1^{i-1}D_{m+1} - (i-2)x_1^iD_{2m+1}, \quad D_K(x_{2m+1}) =  \sum_{j=1}^{2m} x_jD_j + 2x_{2m+1}D_{2m+1}, \]
    and so 
    \[ [D_K(x_{2m+1}), D_K(x_1^i)] = i(i-1)x_1^{i-1}D_{m+1} - i(i-2)x_1^iD_{2m+1} \in \mb{K}_n.\]
    By taking linear combinations, it follows that $x_1^iD_{m+1} \in \mb{K}_n$ for all $i \geq 2$, and all such elements commute. Now,
    \[ D_K(x_1x_{2m+1}) = x_1^2D_1 + x_1 \sum_{j=2}^{2m} x_jD_j + x_{2m+1}D_{m+1} + x_1x_{2m+1}D_{2m+1}, \]
    so let $X = [D_K(x_{2m+1}), D_K(x_1x_{2m+1})]$, which has the expression
    \[ X = x_1^2D_1 + x_1\sum_{j=2}^{2m} x_jD_j + x_{2m+1}D_{m+1}.  \]
    Then $X$ is a homogeneous element of $\mb{K}_n$ (since $x_{2m+1}$ has degree 2) and we can compute that
    \[ [X, x_1^iD_{m+1}] = ix_1^{i+1}D_{m+1}.\]
    Thus, where $Y_i = x_1^iD_{m+1}$, the subalgebra of $\mb{K}_n$ generated by $\{X, Y_i \mid i \geq 2\}$ is graded, and isomorphic to the shift algebra by \Cref{recognise shift lemma}.
\end{proof}

\begin{proof}[Proof of \Cref{graded Lie algebras non-coherence thm}]
Let $\mf{g}$ be an infinite-dimensional $\mb{Z}$-graded simple Lie algebra with polynomially bounded growth. Then, by \Cref{graded Lie algebras classification prop}, $\mf{g}$ is either a loop algebra, of Cartan type, or the Witt algebra.

If $\mf{g}$ is an (untwisted) loop algebra then $\mf{g} = L(\mf{h})$ with $\mf{h}$ a simple, hence non-nilpotent, finite-dimensional Lie algebra, thus by \Cref{loop coherence thm}, $\Ua(\mf{g})$ is not graded coherent. If $\mf{g}=L(\mf{h})^\sigma$ is a twisted loop algebra, then $\Ua(\mf{g})$ is not graded coherent by \Cref{twisted loop algebra non-coherence thm}.

If $\mf{g}$ is of Cartan type, note $\mf{g}$ is not $\Wone$ by hypothesis, so by \Cref{Cartan shift algebra thm}, $\mf{g}$ contains $\mf{s}$ as a graded subalgebra. By \Cref{shift non-coherent prop}, $\Ua(\mf{s})$ is not graded coherent, and thus \Cref{coherence descent lemma} implies $\Ua(\mf{g})$ is not graded coherent. This completes the proof.
\end{proof}

%% file: Graded_coherence_of_B_and_T.tex
\section{Graded (non-) coherence of some rings related to $\Ua(W)$} \label{Graded noncoherence of some rings related to UW}

\addtocounter{onesidedW}{1}

In this section, let $\kk$ be a field of characteristic 0.
As we want to emphasise that the one-sided Witt algebra is a subalgebra of the Witt algebra, we will refer to it as $W_{\geq-1}$ in this section, and sometimes refer likewise to other similar subalgebras, for example $W_{\geq 0}$ and $W_{\geq 1}$. 
We consider several rings related to $\Ua(W)$ or to $\Ua(\Wone)$, and discuss their coherence.

As Lie algebras of longstanding interest, and the only missing cases of \Cref{graded Lie algebras non-coherence thm}, we are motivated to determine the coherence of the universal enveloping algebras of $W$ and $\Wone$. These enveloping algebras were proved to be non-Noetherian in \cite{sierra2014universal}.

In this section we explore the proof strategy of \cite{sierra2014universal} in our new context of coherence. The non-Noetherianity of $\Ua(\Wone)$ is proved via a homomorphism from $\Ua(\Wone)$ to a non-Noetherian idealiser, which is crucially shown to have non-Noetherian image.

We therefore construct a corresponding $\ZZ$-graded and non-Noetherian idealiser $\wt{T}$ and a homomorphism $\Psi: \Ua(W) \to \wt{T}$.
We show that $\Psi(\Ua(W))$ and $\Psi(\Ua(\Wone))$ are not left or right (graded) coherent. Unfortunately, this fails to show coherence or non-coherence of the universal enveloping algebras, due to the kernel of $\Psi$ being an infinitely-generated one-sided ideal.

In contrast, the idealiser $\wt{T}$ and the related ring $T = \wt{T}_{\geq 0}$ are graded right coherent. Though $\wt{T}$ and $T$ are left Noetherian, they are not right Noetherian. We do not know if $\wt{T}$ and $T$ are right coherent, and discuss this problem in \cref{Ungraded coherence and Rees rings section}. Similarly to the above, the graded coherence of $\wt{T}$ is not sufficient to show that $\Ua(W)$ is graded coherent, though it constitutes very weak evidence in this direction. 

\subsection{Non-coherence of $\wtB$ and $B$} \label{Non-coherence of B}

In this subsection we consider the (first) {\em orbit homomorphism } 
\[ \Psi: \Ua(W) \rightarrow \widetilde{A_1} \otimes  \kk[v_0], \quad \Psi(e_n) = t^{n+1}\partial + (n+1)t^nv_0, \]
where $\widetilde{A_1}=\kk[t,t^{-1}][\partial; \frac{d}{dt}]$
is the localised Weyl algebra; 
see \cite{phamXorbit} and \cite{phamXkernel}.
It is known that the image of $\Psi$ is neither left or right Noetherian; see \cite[Proposition 2.5.6]{buzaglo2024lie}, which is the essence of the proof that $\Ua(W)$ is not Noetherian in \cite{sierra2014universal}.

As mentioned in previous sections, it is not known if $\Ua(W)$ is coherent. However, we demonstrate that its image under $\Psi$ is not coherent.

\notinsubfile{ 
\Bnotcoherent* 
}

\onlyinsubfile{
\begin{thm} \label{B not coherent thm}
    The images of $\Ua(W)$ and $\Ua(\Wone)$ under $\Psi$ are not (graded) left or right coherent rings.
\end{thm}
}

Despite \Cref{B not coherent thm}, it is possible that either of the enveloping algebras $\Ua(W)$, $\Ua(\Wone)$ are coherent, since not every quotient of a coherent ring is coherent. It is known that if $R$ is a left coherent ring and $I$ is a two-sided ideal finitely-generated as a left ideal, then $R/I$ is left
coherent, see \cite[Corollary 1.2]{harris1966some1}. However, the kernel of $\Psi$ is not finitely-generated as a left ideal (see Theorem 1.4 of \cite{phamXkernel} for $\Ua(\Wone)$).

We prove \Cref{B not coherent thm} by first giving an explicit description of the image algebras.
Let $\wtB = \Psi(\Ua(W))$ and let $B = \Psi(\Ua(\Wone))$. 
Let $\widetilde S = \widetilde{A_1} \otimes \kk[v_0]$ and let $a = -(t\partial + v_0)$, $b=-v_0$. 
We may write $\widetilde S$ as a skew polynomial ring of automorphism type,
\[ \widetilde S = \kk[a,b][t,t^{-1}; \sigma],\]
where   $\sigma(a)=a+1$, $\sigma(b)=b$. 
(Here our convention  is that if $f \in \kk[a,b]$, then $tf = f^\sigma t$.)
In terms of $a$ and $b$ we have  $\Psi(e_n) = -(a+nb)t^n$ for all $n \in \mb{Z}$. 

If we let $t$ have degree $1$ and elements of $\kk[a,b]$ have degree zero, then $\widetilde S$ is $\mb{Z}$-graded  and
$\Psi$ is a homomorphism of graded $\kk$-algebras, and thus $\wtB, B$ are graded subrings of $\widetilde S$.
We define the homogeneous elements
\[ c_n  = -\Psi(e_n) = (a+nb)t^n, \quad d_n = b(1-b)t^n,\]
and note the following multiplication rules.
\begin{lemma} \label{computations lemma}
Let $m, n \in \mb{Z}$. Then
\begin{flalign*}
&& c_md_n &= (a+mb)d_{m+n},\\
&& d_mc_n &= (a+nb+m)d_{m+n},\\
&& c_mc_n &= (a+m)c_{m+n} - mnd_{m+n},\\
&& d_md_n &=b(1-b)d_{m+n}. && \qed
\end{flalign*}

\end{lemma}

Let us record a particular consequence of \Cref{computations lemma} which we will use more than once.

\begin{lemma} \label{b power d lemma}
    Let $R \in \{\wtB, B\}$. Then $b^{m+1}d_{n} \in R\{b^{m}d_{n}, b^{m}d_{n+1}\}$ for any $m \geq 0$ and $n \in \mb{Z}$.
\end{lemma}
\begin{proof}
    This follows from $c_0d_{n} - c_{-1}d_{n+1} = bd_{n}$, since $b^{m}$ is central and $c_0, c_{-1} \in R$.
\end{proof}

We now compute the graded pieces of the images  $\wtB$ and $B$.

\begin{lemma} \label{B structure lemma}
If $n \in \mb{Z}$ and $n \neq 0$, then
\[ \wtB_n =  \kk[a]c_n+ \kk[a,b]d_n. \]
Further, $\wtB_0 = \kk[a] + \kk[a,b]b(1-b) = \kk + \kk[a]c_0+\kk[a,b]d_0$.
\end{lemma}
\begin{proof}
Let $B'$ be the graded $\kk$-vector subspace of $\wt{S}$ given by $B_n' = \kk[a]c_n + \kk[a,b]d_n$ for $n \neq 0$ and $B_0' = \kk[a]+\kk[a,b]b(1-b)$.

We first show that $\wtB \leq B'$. 
Since $\Ua(W)$ is generated as a $\kk$-algebra by the basis elements $\{e_n \mid n \in \mb{Z}\}$ of $W$, $\wtB$ is generated by $\{ c_n \mid n \in \mb{Z}\}$. 
Clearly $c_n \in B_n'$, so it is then enough to show that $c_nB' \subseteq B'$ for all $n \in \mb{Z}$. Now by \Cref{computations lemma}, for all $n,m \in \mb{Z}$, we have
\begin{align*} c_n\kk[a]c_m= \kk[a]c_nc_m = \kk[a]((a+n)c_{n+m} - mnd_{n+m}) \subseteq B_{n+m}',
\end{align*}
and
\[ c_n\kk[a,b]d_m = \kk[a,b]c_nd_m =\kk[a,b](a+nb)d_{n+m} \subseteq B_{n+m}'. \]
Also $c_n \kk = \kk c_n \subseteq B_n'$. Therefore $c_nB' \subseteq B'$, as required, so $\wtB \leq B'$.

We now show that $B' \leq \wtB$. Notice that $\Psi(e_0)=-a \in \wtB$, therefore $\kk[a]\Psi(e_n) = \kk[a]c_n \subseteq \wtB$ for all $n \in \mb{Z}$. 
By Lemma 2.5.5 of \cite{buzaglo2024lie}, $b(1-b)\widetilde S \subseteq \wtB$, therefore $\kk[a,b]d_n \subseteq \wtB$ for all $n$. 
Since the base field $\kk$ is also contained in $\wtB$, we have $B' \leq \wtB$ as required. 
\end{proof}

\begin{lemma} \label{B-1 structure lemma}
    If $n \geq 0$, then $B_n=\wtB_n$ whilst if $n > 0$, then $B_{-n}=\wtB_0c_{-1}^{n}$.
\end{lemma}
\begin{proof}
We first   show that $B_{n} = \wtB_n$ for all $n \geq 0$. By \Cref{computations lemma},
    \[ (c_0-1)c_n - c_{-1}c_{n+1} = (n+1)d_n,\]
    thus $d_n \in B$ for all $n \geq 0$, as $\kk$ is of characteristic zero. By \Cref{b power d lemma} and induction, $b^md_n \in B$ for all $m, n \geq 0$. Since $a =c_0 \in B$, it follows that
    \[ \wtB_n = \kk[a]c_n + \kk[a][b]d_n \subseteq B_n,\]
    for all $n \geq 0$, giving the claim.
    
    The universal enveloping algebra $\Ua(\Wone)$ is generated by $e_{-1}$ and elements of non-negative degree, meaning that $B$ is generated by $c_{-1}$ and $B_n = \wtB_n$ for all $n \geq 0$. Since $B$ is a graded ring and  
    \[ c_{-1}B_0 = c_{-1}(\kk[a] + \kk[a,b]b(1-b)) = B_0c_{-1},\]
    it follows by induction that $B_{-n} = B_0c_{-1}^n$ for all $n \geq 0$.
\end{proof}

\begin{lemma} \label{J generation lemma}
    Let $R \in \{\wtB, B\}$ and let $J = R \cap Rt^{-1}$. Then $J$ is generated by $\{d_n \mid n \in \mb{Z}\} \cap R$ as a left ideal.
\end{lemma}
\begin{proof}
    Let $X = \{d_n \mid n \in \mb{Z}\} \cap R$ and $\II = \{n \in \mb{Z} \mid d_n \in R\}$. It is straightforward to check from \Cref{B structure lemma,B-1 structure lemma} that if $n \in \II$ then $n+1 \in \II$, from which it follows that $J$ contains $X$. 
    Moreover, \Cref{b power d lemma} and induction imply that $b^md_n \in RX$ for all $m \geq 0$ and $n \in \II$. Since $\kk[a]$ is a subring of $R$, it follows that $\kk[a,b]d_n \subseteq RX$ for all $n \in \II$. 
    In fact \Cref{B structure lemma,B-1 structure lemma} imply $\II = \{n \in \mb{Z} \mid R_n=\wtB_n\}$, so if $n \in \II$, then
    \begin{align*} J_nt^{-n} &= \wtB_nt^{-n} \cap \wtB_{n+1}t^{-(n+1)}\\ &= ( \kk[a](a+nb) + \kk[a,b]b(1-b)) \cap ( \kk[a](a+(n+1)b) + \kk[a,b]b(1-b))\\ &= \kk[a,b]b(1-b),\end{align*}
    from which it follows $J_n = \kk[a,b]d_n \subseteq RX$.
    
    If $n \not \in \II$, then $n = -N <0$ and $R = B$. Now, $R_{-N} = \wtB_0c_{-1}^N$ for any $N \geq 0$. Since $\wtB_0c_{-1} = c_{-1}\wtB_0$, it follows
    \[ J_{-N} = c_{-1}\wtB_0c_{-1}^{N-1} \cap c_{-1}\wtB_0c_{-1}^{N-2}t^{-1} = c_{-1}J_{-(N-1)},\]
    so by induction $J_n \subseteq RJ_0 \subseteq RX$, and the result follows.
\end{proof}

We are now almost ready to give the proof of the main result of this subsection. Before doing so, we establish that $\wtB$ and $B$ have symmetric left-right properties.

\begin{lemma} \label{B opposite lemma}
    Let $R \in \{\wtB, B\}$. Then $R$ is isomorphic to its opposite ring (as a graded algebra).
\end{lemma}
\begin{proof}
    Consider the (graded) $\kk$-algebra homomorphism from $\widetilde{S}$ to its opposite ring,
    \[ i: \widetilde{S} \rightarrow \widetilde{S}^{op}, \quad i(a) = -a, \;i(b)=1-b,\; i(t)=t,\; i(t)=t^{-1}, \]
    which is well-defined because
    \[ i(at) = -(a+1)t = i(t)i(a) = i(a)\cdot_{op}i(t)\]
    and $b$ is central in $\widetilde{S}$ and its opposite. Moreover, $i$ is a (graded) isomorphism. Now, $R$ is generated as a $\kk$-algebra by a subset of $\{ c_n \mid n \in \mb{Z} \}$, and
    \[ i(c_n) = i((a+nb)t^n) = i(a+nb) \cdot_{op} i(t^n) = t^n(-a + n(1-b)) = -c_n \]
    therefore $i(R) = R^{op} \subseteq \widetilde{S}^{op}$, which gives the result.
\end{proof}

In particular, $\wtB$ or $B$ is left (graded) coherent if and only if it is right (graded) coherent.

\begin{proof}[Proof of \Cref{B not coherent thm}]
 Let $R \in \{\wtB, B\}$. By \Cref{B opposite lemma}, it is enough to show that $R$ is not (graded) left coherent. Consider the finitely-generated left ideal
 \[ M = Rd_0 + Rd_1. \]
 We have a short exact sequence of $R$-modules
\[ 0 \rightarrow K \rightarrow R^2 \xrightarrow{\phi} M \rightarrow 0,\]
where $\phi(x,y) = xd_0 + yd_1$ and $K = \Ker \phi$. Then
\begin{align*} K &= \{(x,y) \in R^2 \mid xb(1-b) = -yb(1-b)t \}\\ &=  \{(x,y) \in R^2 \mid x = -yt\}\\ & = \{(-yt, y) \mid y \in R \cap Rt^{-1}\}\\ &\cong J.
\end{align*}
We show that $J = R \cap Rt^{-1}$ is not finitely-generated. By \Cref{J generation lemma}, $J=RX$ where
\[X =\{d_n \mid n \in \mb{Z}\} \cap R, \]
so we will show that no finite subset of $X$ generates $J$. If such a subset exists, then in particular, there exists $N \geq 1$ such that
\[ d_n \in \sum_{j=0}^N Rd_j \subseteq \sum_{j=0}^N \wtB d_j  \]
for all $n \geq 0$.
Now, $\wtB$ has a $\kk$-algebra automorphism which sends $c_n \mapsto -c_{-n}$ and $d_n \mapsto d_{-n}$ for all $n$. Thus, we also have that
\[ d_{-n} \in \sum_{j=-N}^0 Rd_j \subseteq \sum_{j=-N}^0 \wtB d_j \]
for all $n \geq 0$, and hence $J$ is generated by a finite subset of $X$ in $\wtB$. 
So, we may assume without loss of generality that $R=\wtB$ and $\II = \ZZ$. 
Suppose there is a finite subset of $X$ which generates $J$. 
Then there exists $N \geq 1$ such that
\[ d_{N+1} \in \sum_{j=-N}^N \wtB_{N+1-j}d_j, \]
using the fact that $\wtB$ is a graded $\kk$-algebra. By \Cref{B structure lemma} and \Cref{computations lemma}, we compute that for all $j \in \{-N, \dots, N\}$,
\begin{align*} \wtB_{N+1-j}d_j &= \kk[a]c_{N+1-j}d_j + \kk[a,b]d_{N+1-j}d_j\\ &= \big(\kk[a](a+(N+1-j)b)+\kk[a,b]b(1-b)\big)d_{N+1}\\ &\subseteq (\kk[a,b]a + \kk[a,b]b)d_{N+1}.\end{align*}
Obviously $d_{N+1} \not \in (\kk[a,b]a + \kk[a,b]b)d_{N+1}$, so we have a contradiction. Thus, no finite subset of $X$ generates $J$ as a left ideal. Since $X$ is a generating set of $J$, it follows that $J$ is not finitely-generated.

Therefore $K$ is not finitely-generated and hence $M$ is not finitely-presented. So $R$ is not left coherent. Since $M$ is additionally a graded left ideal, it follows that $R$ is also not graded left coherent.
\end{proof}

\subsection{Graded coherence of $\wt{T}$ and $T$} \label{Graded coherence of T subsection}

We have seen that the rings $\wtB = \Psi(\Ua(W)$ and $B = \Psi(\Ua(\Wone))$ are neither Noetherian or coherent.
In this section we consider an overring $\wt{T}$ of $\wtB$ and study its coherence.
We also prove the graded coherence of $T = \wt{T}_{\geq 0}$, which is an overring of $\Psi(\Ua(W_{\geq 0}))$.

Let us define the rings of interest in this section.
\begin{defn} \label{T def}
    Recall that  $\wt{S}= \wt{A_1} \otimes \kk[v_0] = \kk[a,b][t, t^{-1}; \sigma]$.
    Let $\wt{I} = \wt{S}a + \wt{S}b$,
     and let $\wt{T} = \{ s \in \wt{S} \mid \wt{I}s \subseteq \wt{I}\}$ be  the idealiser of the left ideal $\wt{I}$ of $\wt{S}$.

     Likewise, let $S = \wt{S}_{\geq 0} = \kk[a,b][t;\sigma]$, let $I = Sa + Sb$, and let $T = \wt{T}_{\geq 0}$ be the idealiser of $I$ in $S$.
\end{defn}

Explicitly, $\wt{T} = \kk \oplus \wt{I}$ and $T = \kk \oplus I$. 
Notice that $\wt{T}$ (respectively, $T$) is a graded ring because the left ideal $\wt{I}$ (respectively, $I$) is graded.

Throughout this subsection, we use parentheses to denote ideals of $\kk[a,b]$, for example $I = S(a,b)$.

It is known that $\wt{T}$ and $T$ are not right Noetherian:  this is proved in a similar way as the proof that $B $ and that $\Ua(\Wone)$ are not right Noetherian and in fact inspired that proof.
In this subsection we prove the following:

\notinsubfile{
\Tgradedcoherent* 
}

\onlyinsubfile{ 
\begin{thm} \label{T graded coherent thm}
    Although $\wt{T} $ and  $T$ are not right Noetherian, both are  left Noetherian and graded right coherent.
\end{thm}
}

However, we do not know if $\wt{T}$ or  $T$ is right coherent (dropping the graded assumption).  
For more on this see \Cref{Ungraded coherence and Rees rings section}.

The following is a straightforward computation: 
\begin{lemma} \label{structure of T lemma}
    The graded pieces of $\wt{T}$ are given by
    \[ \wt{T}_0 = \kk[a,b], \quad \wt{T}_n = t^n(a,b) = (a+n,b)t^n = \wt{I}_n,\]
    for all $n \neq 0$. In particular, $\wt{T}_n$ is a finitely-generated left and right $\wt{T}_0$-module for all $n$. 
    \qed
\end{lemma}

The key step in proving Theorem~\ref{T graded coherent thm} is that there is an easy way to check if a graded right ideal of $\wt{T}$ or $T$ is finitely-generated.

\begin{proposition} \label{fg graded ideals of T prop}
Let $J$ be a graded right ideal of $T$ and let $\wt{J}$ be a graded right ideal of $\wt{T}$.
Then $J$ is finitely-generated if and only if there is an ideal $\FF$ of $ \kk[a,b]$ so that $J_n = \FF(a,b)^{\sigma^n}t^n$ for $n\gg 0$, 
and $\wt{J}$ is finitely-generated if and only if there is an ideal $\FF$ of $ \kk[a,b]$ so that $\wt{J}_n = \FF(a,b)^{\sigma^n}t^n$ for $|n|\gg 0$.
\end{proposition}

\begin{proof}
We give the proof for $\wt{J}$; the proof for $J$ is similar.

$(\Leftarrow)$ Suppose that $\FF$ exists, so there is $L \in \NN$ with $\wt{J}_n = \FF(a,b)^{\sigma^n}t^n = \FF t^n (a,b)$ for $|n| \geq L$.
For $n \geq L+1$ we have
\begin{align*}
\wt{J}_n & \supseteq \FF t^L(a,b) t^{n-L}(a,b) + \FF t^{L+1}(a,b) t^{n-L-1}(a,b) \\
 & = (\FF (a,b)^{\sigma^L} (a,b)^{\sigma^n } + \FF (a,b)^{\sigma^{L+1}} (a,b)^{\sigma^n }) t^n \\
 & = \FF ((a,b)^{\sigma^L}  +  (a,b)^{\sigma^{L+1}} ) (a,b)^{\sigma^n}t^n \\
 & =\FF (a,b)^{\sigma^n}t^n,
 \end{align*}
 since $(a,b)^{\sigma^{L+1}} \neq (a,b)^{\sigma^L}$ are (co-)maximal ideals of $\kk[a,b]$.
 Likewise, if $n \leq -L-1$ then
\[ 
\wt{J}_n \supseteq \FF t^{-L}(a,b) t^{n+L}(a,b) + \FF t^{-L-1}(a,b) t^{n+L+1}(a,b) = \FF(a,b)^{\sigma^n} t^n.
\]
So \[ \wt{J} = (\wt{J}_{-L-1} \oplus \dots \oplus \wt{J}_{L+1})\wt{T}\]
is finitely-generated.

$(\Rightarrow)$
Suppose that $\wt{J}$ is finitely-generated, so there is $M \in \NN$ so that 
\[ \wt{J} = (\wt{J}_{-M} \oplus \dots \oplus \wt{J}_{M})\wt{T}.\]
If $|n| > M$, then
\beq\label{biff} 
\wt{J}_n = \sum_{i=-M}^M \wt{J}_i T_{n-i} = \sum_{i=-M}^M (\wt{J}_i t^{-i}) (a,b)^{\sigma^n} t^n,\eeq
since $n-i \neq 0$ if $|i|\leq M$.
Let $\JJ_n = \wt{J}_n t^{-n}$ for $n \in \ZZ$; note this is an ideal of $\kk[a,b]$.
Then for $|n| > M$, \eqref{biff} can be rewritten
\[ \JJ_n = \sum_{i=-M}^M \JJ_i (a,b)^{\sigma^n} = \bigl(\sum_{i=-M}^M \JJ_i\bigr) (a,b)^{\sigma^n}.\]
So the result holds for  $\FF = \sum_{i=-M}^M \JJ_i$.
\end{proof}

We next give a lemma on how ideals of $\kk[a,b]$ interact with the $\sigma$-orbit of $(a,b)$.

\begin{lemma}\label{lem1.2}
For any ideals $\FF, \GG$ of $ \kk[a,b]$ we have
\beq\label{foo}
\FF (a,b)^{\sigma^n} \cap \GG (a,b)^{\sigma^n} = (\FF \cap \GG) (a,b)^{\sigma^n}
\eeq
for $|n| \gg 0$.
\end{lemma}
\begin{proof}
Let $C = \kk[a,b]$.  For $n \in \ZZ$, let $\mf{q}_n = (a,b)^{\sigma^n}$.

Let $Q_1, \dots, Q_n$ be the  primary components of $\FF$ that are not $(b)$-primary, and write
\[ \FF = (b^k) \cap Q_1 \cap \dots \cap Q_r. \]
Likewise, write 
\[ \GG = (b^\ell) \cap Q'_1 \cap \dots \cap Q'_s,\]
where the $Q'_j$ are primary but not $(b)$-primary.
Without loss of generality, $\ell \geq k$.
%
%
Since  the only irreducible closed subscheme of $\AA^2$ with an infinite intersection with $\ZZ \times 0$ is $V(b)$, there is $N\in\NN$ so that if $|n| \geq N$ then $\mf{q}_n$ is  comaximal with all $Q_i$, $Q'_j$.
Thus, localising, 
\[ (\FF)_{\mf{q}_n} = (b^k)_{\mf{q}_n}, \quad (\GG)_{\mf{q}_n} = (b^\ell)_{\mf{q}_n}\]
for $|n| \geq N$.

Let $|n| \geq N$ and let $\mf{p}$ be a maximal ideal of $C$.
If $\mf{p} = \mf{q}_n$ then 
\[ (\FF (a,b)^{\sigma^n} \cap \GG (a,b)^{\sigma^n})_{\mf{p}}
= (b^k \mf{q}_n)_{\mf{q}_n} \cap (b^\ell \mf{q}_n)_{\mf{q}^n}
= (b^\ell \mf{q}_n)_{\mf{q}^n} = ((\FF \cap \GG) (a,b)^{\sigma^n})_{\mf{p}}.
\]
On the other hand, for $\mf{p} \neq \mf{q}_n$ then
\[(\FF (a,b)^{\sigma^n} \cap \GG (a,b)^{\sigma^n})_{\mf{p}} 
= (\FF\cap \GG)_{\mf{p}} = ((\FF \cap \GG) (a,b)^{\sigma^n})_{\mf{p}}.
\]
Since the stalks of the left hand side and right hand side of \eqref{foo} are equal at every maximal ideal of $C$, we have established the result.
\end{proof}

In \cref{Ideal intersections in a commutative algebra} we prove a similar result in $\kk[a,b,s]$. 

\begin{proof}[Proof of \Cref{T graded coherent thm}]
That $\wt{T}$ is left Noetherian is the opposite-side version of \cite[Theorem~1.4]{ruth}.  
The proof that $T$ is left Noetherian is similar, but we give it here for completeness.
We follow the proof of  \cite[Proposition~5.3]{rogalski2004idealizer}.
Again, let $C = \kk[a,b]$.
Since
\[ S = T + T t = C \oplus \bigoplus_{n\geq 1} t^n\bigl( (a,b) + (a,b)^{\sigma^{-1}} \bigr)\]
is a finitely-generated left $T$-module, by the graded version of 
\cite[Lemma~2.1]{rogalski2004idealizer} it suffices to show that 
$ \Hom_S(S/I, S/J)$
is a Noetherian left $T$-module for any graded left ideal $J$ of $S$.

Let $J$ be a graded left ideal of $S$.
We can identify $\Hom_S(S/I, S/J)_n$ with
\beq\label{eq:Hom}
\faktor{\{ t^n s  \mid s \in C, I t^ns  \subseteq J\}}{J_n}.
\eeq

Let $J$ be generated in degree $\leq N$ and let 
 $\JJ = t^{-N}J_N  \normal C$.
 For $n \geq N$ we have $J_n= t^n\JJ$.

 For all $n \geq N$, the condition $I t^n s  \subseteq J$ is equivalent to 
 $s \in (\JJ: \mf{q}_{-n})$. By increasing $N$ if necessary, we may assume that $\mf{q}_{-n}$ is not an associated prime of $\JJ$ for any $n \geq N$, and thus 
 $(\JJ: \mf{q}_{-n}) = \JJ$.
 This shows that $\Hom_S(S/I, S/J)_{\geq N} = 0$ and it then follows that $\Hom_S(S/I, S/J)$ is a finitely-generated $C$-module  and thus a Noetherian left $T$-module.
 
We now show that $T$ is graded right coherent. 
It  suffices to prove that if $J, K$ are finitely-generated graded right ideals of $T$, then $J\cap K$ is finitely-generated, by \Cref{fg intersection prop}.
By \Cref{fg graded ideals of T prop} there are ideals $\FF, \GG$ of $ C$ so that $J_n = \FF (a,b)^{\sigma^n} t^n$ and $K_n = \GG (a,b)^{\sigma^n} t^n$ for $n \gg 0$.
By \Cref{lem1.2}, for $n \gg 0$ we have
\[ J_n \cap K_n = (\FF\cap \GG)(a,b)^{\sigma^n} t^n.\]
Applying \Cref{fg graded ideals of T prop} to $J\cap K$, the result follows.

The proof that $\wt{T}$ is graded right coherent is similar.
\end{proof}

%% file: Ungraded_coherence_and_Rees_rings.tex
\section{Ungraded coherence and Rees rings} \label{Ungraded coherence and Rees rings section}

As mentioned, we do not know if $T$ or $\wt{T}$ are (ungraded) coherent. 
In this section, we approach the question of ungraded coherence by relating their coherence to graded coherence of Rees rings.  
We prove a general result showing that graded coherence of the Rees ring implies coherence of the filtered ring, \Cref{Rees coherence thm}.

We consider the canonical gradings on the Rees rings of $T$ and $\wt{T}$, but are unable to determine their graded coherence, or otherwise. We show, however, that the Rees rings are graded coherent under the $u$-grading. Whilst the $u$-grading emulates the grading of $T$, $\wt{T}$ in \Cref{Graded noncoherence of some rings related to UW}, establishing $u$-graded coherence requires a non-trivial commutative algebra result on intersections of ideals, see \cref{Ideal intersections in a commutative algebra}. 

In the final \cref{Eventually graded ideals}, we investigate a mild generalisation of finitely generated graded ideals of $T$, eventually graded ideals, and show they are finitely-presented. We give two examples demonstrating that this approach fails to determine the coherence of $T$.

\subsection{Coherence of Rees rings} \label{Coherence of Rees rings subsection}

In this subsection, let $A$ be a ring 
equipped with an ascending filtration  $(A_j)_{j \in \ZZ}$.
We assume that the filtration is exhaustive, namely, $A = \bigcup_{j \in \mb{Z}} A_j$. Recall that the Rees ring $\Rees(A)$ of $A$ is the $s$-graded subring of $A[s,s^{-1}]$ with $j$th graded piece $A_js^j$.
If $A_j = 0$ for $j<0$, namely $A$ has an ascending $\NN$-filtration, then $\Rees(A) \subseteq A[s]$.

\notinsubfile{
\Reescoherence*
}

\onlyinsubfile{
\begin{thm} \label{Rees coherence thm}
If the Rees ring of $A$ is graded left coherent, then $A$ is left coherent.
\end{thm}
}

\begin{proof}
Throughout this proof, denote by $\phi$ the ring homomorphism $\Rees(A) \rightarrow A$ which takes $as^m$ to $a$. Since the filtration on $A$ is exhaustive, $\phi$ is surjective. In fact, every element of $A$ is the image of a homogeneous element under $\phi$.

Let $I \subseteq A$ be a left ideal of $A$ generated by finitely many elements $S \subseteq A$, and let $f: A^{|S|} \rightarrow I$ be the corresponding left $A$-module homomorphism.
Since the filtration on $A$ is exhaustive, for each $a \in S$ there exists $d_a \in \mb{Z}$ with $as^{d_a} \in \Rees(A)$. Define the finitely-generated graded left ideal
\[ J = \Rees(A)\{as^{d_a} \mid a \in S\},\]
with corresponding surjection $g: \Rees(A)^{|S|} \rightarrow J$. Clearly $\phi(J)=I$, giving a commutative square with surjective vertical maps,
\[
\begin{tikzcd}
\Rees(A)^{|S|} \arrow[r, "g"] \arrow[d, "\phi^S"] & J \arrow[d, "\phi"] \\
A^{|S|} \arrow[r, "f"]                            & I                  
\end{tikzcd} \]
and we let $K = \Ker f$ and $L = \Ker g$. We can compute that
\[ f(\phi^S(L)) = \phi(g(L)) = \phi(0)=0,\]
so $\phi^S(L) \subseteq \Ker f = K$, giving a commutative diagram
\[
\begin{tikzcd}
0 \arrow[r] & L \arrow[r] \arrow[d, "\phi^S|_L"] & \Rees(A)^{|S|} \arrow[r, "g"] \arrow[d, "\phi^S"] & J \arrow[d, "\phi"] \arrow[r] & 0 \\
0 \arrow[r] & K \arrow[r]                        & A^{|S|} \arrow[r, "f"]                            & I \arrow[r]                   & 0.
\end{tikzcd}\]

For any $x \in K$ there exists a homogeneous $y \in \Rees(A)^{|S|}$ with $\phi^S(y)=x$. Now,
\[ \phi(g(y)) = f (\phi^S(y)) = f(x)=0,\]
so $g(y)$ is a graded element of $\Ker \phi$. 
But  the image of any non-zero homogeneous element under $\phi$ is non-zero, from which it follows that $g(y)=0$, and therefore $y \in L$. Thus, $K = \phi^S(L)$.

Now, $\Rees(A)$ is graded left coherent and $J$ is a finitely-generated graded left ideal, so $J$ must be finitely-presented, hence $L$ is finitely-generated. Therefore $K = \phi^S(L)$ is a finitely-generated $A$-module, and $I$ is finitely-presented. It follows that $A$ is left coherent.
\end{proof}

Notice that if $\Rees(A)$ is graded coherent, so is the associated graded ring $\gr A = \Rees(A)/s\Rees(A)$ by \Cref{finitely-presented coherence prop}. However, $\gr A$ being coherent does not imply that $A$ is coherent. If $A$ is any universal enveloping algebra with its usual filtration, $\gr A$ is a polynomial ring, hence is coherent, but $A$ may not be coherent (for example, see \Cref{graded Lie algebras non-coherence thm}).

We mention that if $A$ has a descending filtration, that is $A_0=A$, \Cref{Rees coherence thm} also applies to the descending Rees ring
\[ \Rees(A) \cap A[s^{-1}]  = \bigoplus_{j \leq 0} A_{j}s^{j}  \leq A[s^{-1}],\]
with an identical proof.

\subsection{Coherence results for the Rees ring of $\wt{T}$ and $T$} \label{Coherence results for the Rees ring of T}

Again let $\kk$ be a field of characteristic 0.
As we have said, although $T$ and $\wt{T}$ are graded right coherent, we do not know whether they are right coherent. In view of \Cref{Rees coherence thm}, we now consider the graded coherence of the Rees rings of $T$ and $\wt{T}$.

\begin{defn} \label{Rees T def}
    We define $\mc{R}$ to be the Rees ring $\Rees(T)$ of $T$, where $T$ has the natural increasing filtration induced by its $t$-grading (described in \cref{Non-coherence of B}). That is, $T$ has filtered pieces
    \[ F^n(T) = \bigoplus_{i \leq n} T_i, \]
    so that
    \[ \mc{R} = \Rees(T) = \bigoplus_{n \geq 0} F^n(T)s^n \subseteq T[s].\]
    Similarly, we define $\wt{\mc{R}}$ to be the Rees ring $\Rees(\wt{T})$.
\end{defn}

Of course, $\mc{R}$ and $\wt{\mc{R}}$ are graded with the natural Rees ring grading, and $\mc{R}$ is a graded subring of $\wt{\mc{R}}$. In fact, since $T = \wt{T}_{\geq 0}$, we have that $\mc{R} = \wt{\mc{R}}_{\geq 0}$.

Whilst we do not know if either Rees ring is graded right coherent with the natural grading, we show in the remainder of this section that $\mc{R}$ and $\wt{\mc{R}}$ are graded right coherent with respect to an alternative grading.

Recall from \Cref{structure of T lemma} the graded pieces of $T$ and $\wt{T}$, so that \Cref{Rees T def} gives

\begin{align*} \wt{\mc{R}} &= \kk[a,b,s] \oplus \bigoplus_{\substack{n \neq 0\\m \geq n}} (a+n, b)t^ns^m\\ &= \kk[a,b,s] \oplus  \bigoplus_{\substack{n \neq 0\\i \geq 0}} (a+n, b)u^n s^i, \end{align*}
and so
\begin{align*} \mc{R} = \wt{\mc{R}}_{\geq 0}&= \kk[a,b,s] \oplus  \bigoplus_{\substack{n > 0\\i \geq 0}} (a+n, b)u^n s^i, \end{align*}
where $u=ts$.

\begin{defn}
    The $u$-grading on $\wt{\mc{R}}$ has graded pieces
    \[ \wt{\mc{R}}_n = \begin{cases} \kk[a,b,s] & n=0 \\ \kk[s](a+n, b)u^n & n \neq 0. \end{cases},\]
    and the $u$-grading on $\mc{R}$ is induced from $\wt{\mc{R}}$.
\end{defn}

It is straightforward to check that this makes $\wt{\mc{R}}$ a $\mb{Z}$-graded ring, and $\mc{R}$ an $\NN$-graded ring. We will show this makes $\mc{R}, \wt{\mc{R}}$ graded coherent.

\notinsubfile{ 
\ugradedcoherence*
}

\onlyinsubfile{
\begin{thm} \label{u-graded coherence thm}
    The ring $\mc{R}$ is graded coherent with respect to the $u$-grading on $\mc{R}$.
\end{thm}
}

Note that $\wt{\mc{R}}$ also has a bi-grading which extends the above grading,
\[\wt{\mc{R}}_{i,n} = \begin{cases} \kk[a,b]s^i & n=0 \\ (a+n,b)u^ns^i & n \neq 0, \end{cases} \]
and it is immediate that $\mc{R}$ is also bi-graded coherent.

\begin{corollary} \label{bi-graded coherence cor}
The bi-graded rings $\mc{R}$ and $\wt{\mc{R}}$ are bi-graded right coherent. That is, every finitely-generated bi-graded right ideal of $\mc{R}$ or $\wt{\mc{R}}$ is a finitely-presented right $\mc{R}$ or $\wt{\mc{R}}$-module.
\end{corollary}

\subsection{Ideal intersections in a commutative algebra} \label{Ideal intersections in a commutative algebra}

Recall that the proof that $T$ is graded right coherent, \Cref{T graded coherent thm}, requires the control of intersections of certain ideals in the polynomial ring $\kk[a,b]$, see \Cref{lem1.2}. In this subsection, we prove \Cref{eventual intersection thm}, an analogue of \Cref{lem1.2} for the commutative polynomial ring $\wh{C}= \kk[a,b,s]$. 
This result is  crucial to the proof of \Cref{u-graded coherence thm}.

Again let $\kk$ be a field of characteristic zero.
Let $\sigma$ be the $\kk$-algebra automorphism which is the identity on $\kk[b,s]$ and sends $a \mapsto a+1$. Then $\sigma$ acts on the ideals of $\wh{C}$. Let $I_0 = (a,b)$ and define $I_n = \sigma^n(I_0) = (a+n, b)$ for all $n \geq 0$.

Given two ideals $J,K \subseteq \wh{C}$, we are interested in when the equality
\[ JI_n \cap KI_n = (J \cap K)I_n \]
holds. 
The following lemma allows us to describe this in terms of Tor groups.

\begin{lemma} \label{Tor cokernel lemma}
Let $B$ be a ring and $L$ be a two-sided ideal of $B$. Let $F$ be a right $B$-module and $J, K \leq F$ be submodules. Then the $B$-module
\[ X = \faktor{JL \cap KL}{ (J \cap K)L} \]
is isomorphic to the cokernel of the natural map
\[ \Tor_1^B(J, B/L) \oplus \Tor_1^B(K, B/L) \rightarrow \Tor_1^B(J+K, B/L) \]
induced by the addition map $J \oplus K \xrightarrow{+} J+K$.

\end{lemma}
\begin{proof}
Consider the short exact sequence
\[ 0  \rightarrow J \cap K \xrightarrow{i} J \oplus K \xrightarrow{+} J+K \rightarrow 0, \]
where $i(x) = (x,-x)$, and apply the functor $\otimes_B (B/L)$. We obtain an exact sequence
\begin{equation} \label{E1}\begin{split} 
&\Tor_1^B(J, B/L) \oplus \Tor_1^B(K, B/L) \xrightarrow{g} \Tor_1^B(J+K, B/L) \xrightarrow{h}\\ \xrightarrow{h} &\faktor{J \cap K}{(J \cap K)L} \xrightarrow{i'} \faktor{J}{JL} \oplus \faktor{K}{KL}.
\end{split} \end{equation}

Note that $g$ is the natural map discussed in the lemma statement.
Then,
\[ \Ker i' = \faktor{JL \cap KL \cap (J \cap K)}{(J \cap K)L} = \faktor{JL \cap KL}{(J \cap K)L} = X. \]
From the exactness of \eqref{E1}, we have $X \cong \coker g$.
\end{proof}

Returning to our specific scenario, we see that computing the Tor group $\Tor_1^{\wh{C}}(A, \wh{C}/I_n)$ for ideals $A$ of $\wh{C}$ is of interest. We show that this Tor group embeds in another.

\begin{proposition} \label{injective hom prop}
   Let $n \geq 0$ and let $F$ be a flat $\wh{C}$-module. For any submodule $A \leq F$, there is an injective $\wh{C}$-module homomorphism
   \[ \Tor_1^{\wh{C}}(A, \wh{C}/I_n) \rightarrow \Tor_1^{\wh{C}}(F/A, \wh{C}/(a+n)). \]
\end{proposition}
\begin{proof}
There is an exact sequence of $\wh{C}$-modules
\beq\label{Edag} 0 \rightarrow (a+n) \rightarrow I_n \rightarrow I_n/(a+n) \rightarrow 0. \eeq
Since $a+n \in \wh{C}$ is not a zero divisor, $(a+n)$ is a free, hence flat, $\wh{C}$-module. Therefore, upon applying the functor $(F/A) \otimes_{\wh{C}}$ to \eqref{Edag}, we obtain an exact sequence which begins
\beq\label{Eddag} 0 \rightarrow \Tor_1^{\wh{C}}(F/A, I_n) \rightarrow \Tor_1^{\wh{C}}(\wh{F}/A, I_n/(a+n)). \eeq
Now, 
\[ I_n/(a+n) = (a+n, b)/(a+n) \cong (b)/(b(a+n)) \cong \wh{C}/(a+n), \]
whilst dimension shifting shows that
\[  \Tor_1^{\wh{C}}(F/A, I_n) \cong  \Tor_2^{\wh{C}}(F/A, \wh{C}/I_n) \cong  \Tor_1^{\wh{C}}(A, \wh{C}/I_n), \]
because $F$ is flat. Thus \eqref{Eddag} gives an injective homomorphism as in the statement. 
\end{proof}

\begin{corollary} \label{injective hom ideal cor}
Let $n \geq 0$. For any ideal $A \subseteq \wh{C}$, there is an injective $\wh{C}$-module homomorphism
\[ \Tor_1^{\wh{C}}(A, \wh{C}/I_n) \rightarrow \Tor_1^{\wh{C}}(\wh{C}/A, \wh{C}/(a+n)). \]
\qed
\end{corollary}

We now show that the second Tor group can be non-zero only finitely many times.

\begin{proposition} \label{finitely many prop}
Let $M$ be a Noetherian $\wh{C}$-module. There are only finitely many  integers $n$ such that $\Tor_1^{\wh{C}}(M, \wh{C}/(a+n))$ is non-zero.
\end{proposition}
\begin{proof}
For every $n$, because $a+n$ is not a zero divisor of $\wh{C}$, one can compute that
\begin{align*} \Tor_1^{\wh{C}}(M, \wh{C}/(a+n)) &\cong  \{m \in M \mid (a+n)m=0 \} =:Y_n. \end{align*}
Notice that if $m \in Y_n \cap Y_p$ for $n \neq p$, then
\[ (n-p)m = (a+n)m - (a+p)m = 0, \]
so $m=0$ because $\kk$ is of characteristic zero. Hence $Y_n \cap Y_p = 0$.

Now let $n, p_1, \dots, p_N$ be distinct non-negative integers. 
Let
\[ m \in Y_n \cap (Y_{p_1} + \dots + Y_{p_N}), \]
and write $m=m_1 + \dots + m_N$ where all $m_j \in Y_j$.
Then
\[ \sum_{j=1}^N (n-p_j)m_j = \sum_{j=1}^N (a+p_j)m_j + (n-p_j)m_j = (a+n)m=0.\]
Because $n \neq p_N$ and $\kk$ is of characteristic zero, it follows that $m_N \in Y_{p_1} + \dots + Y_{p_{N-1}}$. So, 
\[ m \in Y_n \cap (Y_{p_1} + \dots + Y_{p_{N-1}}),\]
and thus by induction $m=0$. Therefore,
\[ Y_n \cap (Y_{p_1} + \dots + Y_{p_N})=0. \]
Since this holds for any choice of distinct $n, p_1, \dots, p_N$, it follows that the submodule of $M$ generated by the $Y_n$ is a direct sum,
\[ \sum_{n \geq 0} Y_n = \bigoplus_{n \geq 0} Y_n. \]
Since $M$ is Noetherian, $Y_n=0$ for all but finitely many $n$. Thus $\Tor_1^{\wh{C}}(M, \wh{C}/(a+n))=0$ for all but finitely many $n$.
\end{proof}

\begin{corollary} \label{finitely many ideal cor}
Let $A \subseteq \wh{C}$ be any ideal. There are only finitely many integers $n$ such that $\Tor_1^{\wh{C}}(\wh{C}/A, \wh{C}/(a+n))$ is non-zero. \qed
\end{corollary}

\begin{corollary} \label{Tor eventually zero cor}
Let $F$ be a Noetherian flat $\wh{C}$-module and $A \leq F$ be any submodule. For all but finitely many $n \in \mb{Z}$, $\Tor_1^{\wh{C}}(A, \wh{C}/I_n)=0$.
\end{corollary}
\begin{proof}
The module $\Tor_1^{\wh{C}}(A, \wh{C}/I_n)$ is isomorphic to a submodule of $\Tor_1^{\wh{C}}(F/A, \wh{C}/(a+n))$ for all $n$, by \Cref{injective hom prop}. 
\Cref{finitely many prop} then implies $\Tor_1^{\wh{C}}(F/A, \wh{C}/(a+n))=0$ for all but finitely many $n$, forcing $\Tor_1^{\wh{C}}(A, \wh{C}/I_n)=0$.
\end{proof}

\begin{corollary} \label{Tor eventually zero ideal cor}
Let $A \subseteq \wh{C}$ be any ideal. For all but finitely many $n \in \mb{Z}$, the Tor group $\Tor_1^{\wh{C}}(A, \wh{C}/I_n)=0$. \qed
\end{corollary}

\Cref{Tor eventually zero cor} has a further corollary. 
Recall that $\mc{R} = \Rees(T)$.  
For any $n \geq 0$, let $G_n \leq \mc{R}$ be the right $\wh{C}$-module given by the first $n+1$ $u$-graded pieces of $\mc{R}$, namely,
\beq\label{G} G_n = \bigoplus_{i=0}^n (u^i\wh{C} \cap \mc{R}).\eeq
Now, $u^i\wh{C} \cap \mc{R} = u^i(a,b)$ unless $i=0$, hence
\[ G_n \cong \wh{C} \oplus (a,b)^{\oplus n}.\]

\begin{corollary} \label{Tor eventually zero u grading cor}
    Let $m \geq 0$ and $A \leq G_m$. For all but finitely many integer $n$, we have $\Tor_1^{\wh{C}}(A, \wh{C}/I_n)=0$.
\end{corollary}
\begin{proof}
    Since $G_m \cong \wh{C} \oplus (a,b)^{\oplus m} \leq \wh{C}^{m+1}$, we may view $A$ as a submodule of $\wh{C}^{m+1}$. Since $\wh{C}^{m+1}$ is a Noetherian flat $\wh{C}$-module, we can apply \Cref{Tor eventually zero cor} to obtain the result.
\end{proof}

We can now prove the main result of this subsection.

\begin{thm} \label{eventual intersection thm}
Let $F$ be a Noetherian flat $\wh{C}$-module and $J, K \leq F$. For all but finitely many integer $n$,
\[ JI_n \cap KI_n = (J \cap K)I_n. \]
\end{thm}
\begin{proof}
Let $n \geq 0$. By \Cref{Tor cokernel lemma}, $JI_n \cap KI_n = (J \cap K)I_n$ if the natural map
\[ \Tor_1^{\wh{C}}(J, \wh{C}/I_n) \oplus \Tor_1^{\wh{C}}(K, \wh{C}/I_n) \rightarrow \Tor_1^{\wh{C}}(J+K, \wh{C}/I_n) \]
is surjective. By \Cref{Tor eventually zero cor}, for all but finitely many $n$, the Tor group $\Tor_1^{\wh{C}}(J+K, \wh{C}/I_n)$ is zero, forcing the map to be surjective. This gives the result. 
\end{proof}

\begin{corollary} \label{eventual intersection ideal cor}
Let $J,K \subseteq \wh{C}$ be ideals. For all but finitely many integer $n$,
\[ JI_n \cap KI_n = (J \cap K)I_n. \]
\end{corollary}
\begin{proof}
Notice again that $\wh{C}$ is a Noetherian flat $\wh{C}$-module and apply \Cref{eventual intersection thm}.
\end{proof}

Finally, recall the modules $G_n$ from \eqref{G}.

\begin{corollary}
Let $m \geq 0$ and $J, K \leq G_m$. For all but finitely many integer $n$,
\[ JI_n \cap KI_n = (J \cap K)I_n.\]
\end{corollary}
\begin{proof}
Since $G_m \cong \wh{C} \oplus (a,b)^{\oplus m} \leq \wh{C}^{m+1}$, we may view $J,K$ as submodules of $\wh{C}^{m+1}$. 
We then apply \Cref{eventual intersection thm} to obtain the result.
\end{proof}

Notice that the Tor groups $\Tor_1^{\wh{C}}(M, \wh{C}/I_n)$, for $M$ a $\wh{C}$-module, play a crucial role in the arguments of this section, and therefore in the graded coherence of the Rees ring of the idealiser $T$. We comment that this is reminiscent of the role that Tor groups play in the (non-)Noetherianity of idealisers, as discussed in \cite{S-idealizer} and \cite{ruth}.

\subsection{Proof of \Cref{u-graded coherence thm}} \label{Proof of u-graded coherence thm subsection}

We now prove that the Rees algebras $\mc{R} = \Rees(T)$ and $\wt{\mc{R}} = \Rees(\wt{T})$ are $u$-graded coherent.

\begin{proof}[Proof of \Cref{u-graded coherence thm}]
    We prove that $\wt{\mc{R}}$ is $u$-graded coherent; the proof for $\mc{R}$ is very similar. The $u$-graded pieces of $\wt{\mc{R}}$ are
    \[ \wt{\mc{R}}_j = \begin{cases} \wh{C}& j=0 \\ \wh{C}(a+j,b)u^j & j \neq 0. \end{cases}\]
    Note that $u^j\wh{C}(a,b) = \wh{C}(a+j,b)u^j = I_ju^j$.
    
    Let $J\leq \wt{\mc{R}}$ be a right ideal which is graded with respect to the $u$-grading. Then, we have the following.
    \begin{claim} \label{fg claim}
       $J$ is a finitely-generated right ideal if and only if, there exists an ideal $\FF \leq \wh{C}$ and $N \in \mb{N}$ such that $J_j = \FF(a+j,b)u^j = \FF I_ju^j$ whenever $|j| \geq N$. 
    \end{claim}
    \begin{proof}
    This can be proved as in the proof of \Cref{fg graded ideals of T prop}, replacing $\wt{T}$ by $\wt{\mc{R}}$, $\kk[a,b]$ by $\wh{C}$, and $t$ by $u$.
    \end{proof}
    Given this, let $J,K$ be two finitely-generated right ideals of $\wt{\mc{R}}$ which are graded under the $u$-grading described above. Then there exist ideals $\FF, \GG \leq \wh{C}$ and $N_J, N_K \in \mb{N}$ as in \Cref{fg claim}. Let $N = \max\{N_J, N_K\}$. Then if $|j| \geq N$, we have
    \[ (J \cap K)_j = J_j \cap K_j = (\FF I_j \cap \GG I_j)u^j. \]
    By \Cref{eventual intersection ideal cor}, there exists $M \in \mb{N}$ such that 
    \[ \FF I_j \cap \GG I_j = (\FF \cap \GG)I_j\]
    if $|j| \geq M$. Thus whenever $|j| \geq \max\{M,N\}$, we have
    \[ (J \cap K)_j = (\FF \cap \GG)I_ju^j, \]
    hence by \Cref{fg claim}, $J \cap K$ is a finitely-generated right ideal of $\wt{\mc{R}}$.
    
    Since the intersection of any two finitely-generated $u$-graded right ideals is finitely-generated, and $\wt{\mc{R}}$ is obviously a domain, it follows from \Cref{fg intersection prop} that $\wt{\mc{R}}$ is graded right coherent under the $u$-grading.
\end{proof}

\subsection{Eventually graded ideals} \label{Eventually graded ideals}

Although the results above are not strong enough to show that $T$ is right coherent, we mention that every finitely-generated \emph{eventually graded} right ideal of $T$ is finitely-presented.

\begin{defn} \label{eventually graded def}
    Let $J \subseteq T$ be a right ideal. We say that $J$ is eventually graded (for the standard $t$-grading) if there exists $N \geq 0$ such that
    \[ J \cap \Big(\bigoplus_{n \geq N} T_n\Big) = \bigoplus_{n \geq N} J \cap T_n.\]
\end{defn}

\begin{lemma} \label{finite presentation transfer lemma}
    Let $R, S$ be rings, let $\phi: R \rightarrow S$ be a homomorphism such that $S$ is a finitely-presented right $R$-module, and let $M$ be a right $S$-module. If $M$ is finitely-presented as an $S$-module, then $M$ is finitely-presented as an $R$-module.
\end{lemma}
\begin{proof}
    Since $M$ is finitely-presented, there exists a short exact sequence of $S$-modules
    \[ S^{n} \rightarrow S^m \rightarrow M \rightarrow 0.\]
    Any finitely-generated free $S$-module is finitely-presented as an $R$-module, by assumption, so it follows from Theorem 2.1.2(2) of \cite{glaz1989commutative} that $M$ is a finitely-presented $R$-module.
\end{proof}

\begin{proposition}
    If $J$ is a finitely-generated eventually graded right ideal of $T$, then $J$ is finitely-presented.
\end{proposition}
\begin{proof}
    Let $J$ be a finitely-generated eventually graded right ideal of $T$. Let $N \geq 0$ be as in \Cref{eventually graded def}, and let 
    \[K_N = \bigoplus_{n \geq N} T_n\]
    be the corresponding two-sided ideal, so that
    \[ J' = J \cap K_N = \bigoplus_{n \geq N} J \cap T_n.\]
    Then $J'$ is a graded right ideal of $T$, contained in $J$. Also, $JK_N \subseteq  J \cap K_N \subseteq J'$, from which it follows that $J/J'$ is a right $T/K_N$-module.
    
    Now, $T_0 = \kk[a,b]$ is a Noetherian ring, and $T_n$ is a finitely-generated right $T_0$-module for all $n \geq 0$ by \Cref{structure of T lemma}, from which it follows that $T/K_N$ is a Noetherian right $T_0$-module, and thus a right Noetherian ring.
    
    Since $J$ is a finitely-generated $T$-module, $J/J'$ is a finitely-generated $T/K_N$-module. Thus $J/J'$ is a finitely-presented right $T$-module because $T/K_N$ is Noetherian and by \Cref{finite presentation transfer lemma}. Since $J$ is finitely-generated, it follows that $J'$ is also finitely-generated by \cite[Theorem 2.1.2(3)]{glaz1989commutative}. But $J'$ is graded, thus by \Cref{T graded coherent thm}, $J'$ is finitely-presented. Since $J'$ and $J/J'$ are both finitely-presented, $J$ must be finitely-presented. 
\end{proof}

\begin{example}
    Let $J = (ta+b)T + tbT$. Then $J$ is a finitely-generated right ideal of $T$. We can compute that for all $m \geq 1$ and $f \in (a,b)$,
    \[ tat^mf = (ta+b)t^mf -tbt^{m-1}f \in J,\]
    so $taT_{\geq 1} \subseteq J$, and 
    \begin{align*} J &= (ta+b)T + tbT\\ &= (ta +b)T_0 + (ta+b)T_{\geq 1} + tbT\\ &= (ta +b)\kk[a,b] + taT_{\geq 1} + tbT.\end{align*}
    Thus if $x \in J \cap T_{\geq 2}$, we must have $x \in taT_{\geq 1} + tbT$, and thus is a sum of homogeneous elements of $J$. Hence, $J$ is eventually graded, but is not a graded ideal because $ta+b$ is not a sum of homogeneous elements of $J$.
\end{example}

Any principal right ideal generated by a non-homogeneous element is not eventually graded. However, such ideals are trivially finitely-presented since $T$ is a domain. There also appear to be numerous examples of non-principal right ideals which are not eventually graded. 

\begin{example}
    Let $J = (ta+a)T + tbT$. Then, $J$ is a finitely-generated right ideal which is not eventually graded. To see this, notice that for any $m \geq 0$ and $f \in (a,b)$, we have
    \[ (a+ta)t^mf = t^m(a-m)f + t^{m+1}(a-m)f \in J.\]
    If $J$ were eventually graded, then it would follow that for all large enough $m$ and all $f \in (a,b)$,
    \[ t^m(a-m)f - t^m(a-(m-1))f = t^mf \in J.\]
    However, we can check that $t^ma \not \in J$ for all $m \geq 0$, hence $J$ is not eventually graded.
    
    To see this last statement, we can compute that a general element $z \in J$ is of the form
    \begin{align}\begin{split} \label{ideal general element eq} z&= (ta+a)\sum_{i \geq 0}t^ix_i + tb\sum_{i \geq 0}t^iy_i\\ &=  x_0a + \sum_{i \geq 1} t^i( (a-i)x_i + (a-(i-1))x_{i-1} + by_{i-1} )\end{split}\end{align}
    for some $x_0, y_0 \in \kk[a,b]$, $x_i, y_i \in (a,b)$ for $i \geq 1$. If $a$ or $ta$ lies in $J$, then $ a \in J$, but putting $z=a$ above forces
    \[ (a-1)x_1 + a + by_0=0,\]
    which cannot occur, by considering the equation modulo $(b)$. If $z= t^ma \in J$ for some $m \geq 2$, considering the equation for $z$ modulo $(b)$ shows that $x_i \in (b)$ for all $i < m$, but then we obtain the equation
    \[ (a-m)x_m=a\]
    in $\kk[a]$, which is impossible.

    We now determine whether $J$ is finitely-presented. Putting $z=0$ in \eqref{ideal general element eq}, we quickly find  that $x_0=0$, for all $i \geq 1$ there exists $w_i \in \kk[a,b]$ with $x_i=bw_i$, and that
    \[ (a-(i+1))w_{i+1} + (a-i)w_i \in (a,b) \]
    for all $i \geq 1$. Since $w_i=0$ for all large enough $i$, we deduce that $w_i \in (a,b)$ for all $i \geq 1$, and hence
    \[ \sum_{i \geq 0}t^ix_i \in tb(Sa+Sb) = tbJ.\]
    Conversely,
    \[ (ta+a)tbJ = tb(t+1)(a-1)(Sa+Sb) \subseteq tb T,\]
    therefore the right $T$-module
    \[ \{f \in T \mid (ta+a)f \in tbT \} = tbI.\]
    Since $T$ is a domain, it follows that the relation module for the obvious presentation of $J$ is isomorphic to $I$. Since $I$ is a finitely-generated right ideal, generated by $\{a,b,ta,tb\}$, it follows that $J$ is finitely-presented.

    Thus $J$ is not a counterexample to our expectation that $T$ is right coherent.
\end{example}

%% file: main.bib
@Article{phamXkernel,
  author           = {Tuan Anh Pham and James Timmins},
  journal          = {arXiv preprint},
  title            = {{The Kernel and Image of Orbit Homomorphisms for the Witt Algebra}},
  note             = {\href{https://arxiv.org/abs/2510.00756v1}{arXiv:2510.00756v1}},
  creationdate     = {2026-01-13T16:50:29},
  modificationdate = {2026-01-13T16:50:29},
}

@Article{phamXorbit,
  author           = {Tuan Anh Pham},
  journal          = {arXiv preprint},
  title            = {{The orbit method for the Virasoro algebra}},
  note             = {\href{https://arxiv.org/abs/2504.14670v1}{arXiv:2504.14670v1}},
  creationdate     = {2026-01-13T16:59:50},
  modificationdate = {2026-01-13T16:59:50},
}

@PhdThesis{buzaglo2024lie,
  author  = {Lucas Buzaglo},
  school  = {University of Edinburgh},
  title   = {Lie algebras of derivations and their universal enveloping algebras},
  year    = {2024},
  note    = {\href{http://dx.doi.org/10.7488/era/4901}{DOI:10.7488/ERA/4901}},
  comment = {\href{http://dx.doi.org/10.7488/era/4901}{DOI:10.7488/ERA/4901}},
  doi     = {10.7488/ERA/4901},
}

@Article{sierra2014universal,
  author           = {Susan J. Sierra and Chelsea Walton},
  journal          = {Advances in Mathematics},
  title            = {{The universal enveloping algebra of the Witt algebra is not noetherian}},
  year             = {2014},
  pages            = {239-260},
  volume           = {262},
  creationdate     = {2026-01-20T16:41:51},
  modificationdate = {2026-01-20T16:43:34},
}

@Article{mathieu1992classification,
  author           = {Olivier Mathieu},
  journal          = {Inventiones mathematicae},
  title            = {{Classification of simple graded Lie algebras of finite growth}},
  year             = {1992},
  pages            = {455-519},
  volume           = {108},
  creationdate     = {2026-01-27T16:26:17},
  modificationdate = {2026-01-27T16:27:21},
}

@Article{fialowski1983classification,
  author           = {A. Fialowski},
  journal          = {Moscow University Mathematics Bulletin},
  title            = {{Classification of Graded Lie Algebras with Two Generators}},
  year             = {1983},
  number           = {2},
  pages            = {76-79},
  volume           = {38},
  creationdate     = {2026-02-17T12:08:02},
  modificationdate = {2026-02-17T13:41:13},
}

@Article{dusautoy2004motivic,
  author           = {{du Sautoy}, Marcus and Loeser, François},
  journal          = {Selecta Mathematica},
  title            = {{Motivic zeta functions of infinite-dimensional Lie algebras}},
  year             = {2004},
  number           = {253},
  volume           = {10},
  creationdate     = {2026-02-17T15:31:05},
  doi              = {10.1007/s00029-004-0361-y},
  modificationdate = {2026-02-17T15:37:16},
}

@Article{harris1966some1,
  author           = {Morton E. Harris},
  journal          = {Proceedings of the American Mathematical Society},
  title            = {{Some Results on Coherent Rings}},
  year             = {1966},
  pages            = {474-479},
  volume           = {17},
  modificationdate = {2025-11-21T11:31:55},
}

@Book{glaz1989commutative,
  title     = {Commutative Coherent Rings},
  publisher = {Springer-Verlag},
  year      = {1989},
  author    = {Sarah Glaz},
  number    = {1371},
  series    = {Lecture Notes in Mathematics},
}

@Article{rogalski2004idealizer,
  author           = {Daniel Rogalski},
  journal          = {Journal of Algebra},
  title            = {Idealizer rings and noncommutative projective geometry},
  year             = {2004},
  pages            = {791-809},
  volume           = {279},
  creationdate     = {2026-03-26T13:25:31},
  doi              = {10.1016/j.jalgebra.2004.04.001},
  modificationdate = {2026-03-26T13:26:46},
}

@Article{biswal2025ideals,
  author           = {Rekha Biswal and Susan J. Sierra},
  journal          = {Algebra \& Number Theory},
  title            = {{Ideals in enveloping algebras of affine Kac–Moody algebras}},
  year             = {2025},
  number           = {6},
  pages            = {1199-1230},
  volume           = {19},
  creationdate     = {2026-06-02T13:00:20},
  doi              = {https://doi.org/10.2140/ant.2025.19.1199},
  modificationdate = {2026-06-02T13:01:23},
}

@Book{amayo1974infinite,
  author           = {Ralph K. Amayo and Ian Stewart},
  publisher        = {Noordhoff International Publishing, Leyden},
  title            = {{Infinite-dimensional Lie algebras}},
  year             = {1974},
  creationdate     = {2025-09-26T15:00:19},
  modificationdate = {2025-09-26T15:01:54},
}

@Article{andruskiewitsch2025noetherian,
  author           = {Nicolás Andruskiewitsch and Olivier Mathieu},
  journal          = {Journal of the Mathematical Society of Japan},
  title            = {{Noetherian enveloping algebras of simple graded Lie algebras}},
  year             = {2025},
  number           = {4},
  pages            = {1233-1247},
  volume           = {77},
  creationdate     = {2025-08-28T10:20:51},
  modificationdate = {2025-09-26T14:57:28},
}

@Article{minamoto2018criterion,
  author           = {Hiroyuki Minamoto},
  journal          = {European Journal of Mathematics},
  title            = {{A criterion for graded coherence of tensor algebras and applications to higher dimensional AR-theory}},
  year             = {2018},
  pages            = {612-621},
  volume           = {4},
  creationdate     = {2026-08-03T11:43:06},
  doi              = {https://doi.org/10.1007/s40879-017-0202-0},
  modificationdate = {2026-08-03T11:44:41},
}

@Article{chase1960direct,
  author  = {Stephen U. Chase},
  journal = {Transactions of the American Mathematical Society},
  title   = {Direct products of modules},
  year    = {1960},
  pages   = {457-473},
  volume  = {97},
}

@Article{sierra2016maps,
  author           = {Susan J. Sierra and Chelsea Walton},
  journal          = {Pacific Journal of Mathematics},
  title            = {Maps from the enveloping algebra of the positive {W}itt algebra to regular algebras},
  year             = {2016},
  number           = {2},
  pages            = {475-509},
  volume           = {284},
  creationdate     = {2025-08-05T16:11:03},
  doi              = {10.2140/pjm.2016.284.475},
  modificationdate = {2025-08-05T16:51:36},
}

@Article{piontkovski2008coherent,
  author           = {Dmitri Piontkovski},
  journal          = {Journal of Algebra},
  title            = {Coherent algebras and noncommutative projective lines},
  year             = {2008},
  number           = {8},
  pages            = {3280-3290},
  volume           = {319},
  creationdate     = {2025-12-10T13:21:06},
  modificationdate = {2025-12-10T13:21:06},
}

@Article{polishchuk2005noncommutative,
  author           = {A. Polishchuk},
  journal          = {Mathematical Research Letters},
  title            = {{Noncommutative Proj and coherent algebras}},
  year             = {2005},
  number           = {1},
  pages            = {63-74},
  volume           = {12},
  creationdate     = {2025-12-10T13:21:12},
  doi              = {https://dx.doi.org/10.4310/MRL.2005.v12.n1.a7},
  modificationdate = {2025-12-10T13:24:30},
}

@Article{shotton2020category,
  author    = {Shotton, Jack},
  journal   = {Documenta Mathematica},
  title     = {{The category of finitely presented smooth mod $p$ representations of $GL_2(F)$}},
  year      = {2020},
  pages     = {143--157},
  volume    = {25},
}

@Article{timmins2023coherence,
  author    = {Timmins, James},
  title     = {{Coherence of augmented Iwasawa algebras}},
  journal   = {Advances in Mathematics},
  year      = {2023},
  volume    = {417},
  number    = {108916},
}

@InBook{emerton2025categorical,
  author           = {Emerton, Matthew and Gee, Toby and Hellmann, Eugen},
  chapter          = {{An introduction to the categorical $p$-adic {L}anglands program}},
  editor           = {Pierre-Henri Chaudouard and Wee Teck Gan and Tasho Kaletha and Yiannis Sakellaridis},
  pages            = {167-419},
  publisher        = {American Mathematical Society, Providence, RI},
  title            = {{The Langlands Program: II. Geometrization of the Langlands Correspondence}},
  year             = {2025},
  number           = {112.2},
  series           = {Proceedings of Symposia in Pure Mathematics},
  doi              = {https://doi.org/10.1090/pspum/112.2/02066},
  modificationdate = {2026-06-24T11:42:44},
}

@Article{bell2025enveloping,
  author           = {Jason Bell and Lucas Buzaglo},
  journal          = {International Mathematics Research Notices},
  title            = {{Enveloping Algebras of Derivations of Commutative and Noncommutative Algebras}},
  year             = {2025},
  pages            = {1-11},
  volume           = {17},
  creationdate     = {2026-08-04T15:10:44},
  doi              = {https://doi.org/10.1093/imrn/rnaf265},
  modificationdate = {2026-08-04T15:11:47},
}

@article {smith,
    AUTHOR = {Smith, S. Paul},
     TITLE = {Category equivalences involving graded modules over path
              algebras of quivers},
   JOURNAL = {Adv. Math.},
  FJOURNAL = {Advances in Mathematics},
    VOLUME = {230},
      YEAR = {2012},
    NUMBER = {4-6},
     PAGES = {1780--1810},
      ISSN = {0001-8708,1090-2082},
   MRCLASS = {16W50 (16G20)},
  
MRREVIEWER = {Gena\ Puninski},
       DOI = {10.1016/j.aim.2012.03.031},
       URL = {https://doi.org/10.1016/j.aim.2012.03.031},
}

@article {ruth,
    AUTHOR = {Reynolds, Ruth},
     TITLE = {Idealisers in skew group rings},
   JOURNAL = {J. Algebra},
  FJOURNAL = {Journal of Algebra},
    VOLUME = {538},
      YEAR = {2019},
     PAGES = {207--231},
      ISSN = {0021-8693,1090-266X},
   MRCLASS = {16S35 (16S38 16W50)},
  
MRREVIEWER = {Theodora\ Theohari-Apostolidi},
       DOI = {10.1016/j.jalgebra.2019.06.045},
       URL = {https://doi.org/10.1016/j.jalgebra.2019.06.045},
}

@article {S-idealizer,
    AUTHOR = {Sierra, Susan J.},
     TITLE = {Geometric idealizer rings},
   JOURNAL = {Trans. Amer. Math. Soc.},
  FJOURNAL = {Transactions of the American Mathematical Society},
    VOLUME = {363},
      YEAR = {2011},
    NUMBER = {1},
     PAGES = {457--500},
      ISSN = {0002-9947,1088-6850},
   MRCLASS = {16S38},
  
       DOI = {10.1090/S0002-9947-2010-05110-4},
       URL = {https://doi.org/10.1090/S0002-9947-2010-05110-4},
}

@incollection {rogalski,
    AUTHOR = {Rogalski, Daniel},
     TITLE = {Noncommutative projective geometry},
 BOOKTITLE = {Noncommutative algebraic geometry},
    SERIES = {Math. Sci. Res. Inst. Publ.},
    VOLUME = {64},
     PAGES = {13--70},
 PUBLISHER = {Cambridge Univ. Press, New York},
      YEAR = {2016},
      ISBN = {978-1-107-57003-0; 978-1-107-12954-2},
   MRCLASS = {14A22 (16S38)},
  
MRREVIEWER = {Peter\ Goetz},
}
